\documentclass[11pt]{article}

\usepackage{amscd}      
\usepackage{amssymb}     
\usepackage[intlimits]{amsmath}     
\usepackage{xypic}      
\LaTeXdiagrams          
\usepackage[all,v2]{xy}      
\xyoption{2cell}
\UseAllTwocells
\xyoption{frame}
\CompileMatrices

\usepackage{theorem}

\newtheorem{prop}{Proposition}[section]  
\newtheorem{lem}[prop]{Lemma}

\newtheorem{cor}[prop]{Corollary}
\newtheorem{them}[prop]{Theorem}

\theorembodyfont{\upshape}

\newtheorem{defn}[prop]{Definition}

\newtheorem{numrmk}[prop]{Remark}

\newtheorem{numex}[prop]{Example}
\newtheorem{numrmks}[prop]{Remarks}

\newenvironment{pf}{\begin{trivlist}\item[]{\sc Proof.}}%
            {\nolinebreak $\Box$ \end{trivlist}}
\newenvironment{pff}{\begin{trivlist}\item[]{\sc Proof of Proposition \ref{lem:precisions}.}}%
{\nolinebreak $\Box$ \end{trivlist}}
\newenvironment{pf1}{\begin{trivlist}\item[]{\sc Proof of Proposition 
\ref{pro:intert}.}}%
            {\nolinebreak $\Box$ \end{trivlist}}

\newcommand{\noprint}[1]{}

\renewcommand{\tilde}{\widetilde}

\newcommand{\toto}{\rightrightarrows}

\newcommand{\upst}{^{\ast}}

\newcommand{\com}{{\scriptscriptstyle\bullet}}
\newcommand{\lcom}{_{\scriptscriptstyle\bullet}}
\newcommand{\upcom}{^{\scriptscriptstyle\bullet}}

\newcommand{\Gg}{{\mathfrak g}}

\newcommand{\Cc}{{\mathfrak c}}

\newcommand{\zz}{{\mathbb Z}}

\newcommand{\cc}{{\mathbb C}}

\renewcommand{\O}{{\cal O}}

\newcommand{\hH}{{\cal H}}

\newcommand{\del}{\partial}

\newcommand{\pr}{\mathop{\rm pr}\nolimits}

\newcommand{\Ad}{\mathop{\rm Ad}\nolimits}

\newcommand{\id}{\mathop{\rm id}\nolimits}

\newcommand{\ldiag}[1]%
       {\makebox[0cm]{${\scriptstyle#1}\downarrow\phantom{\scriptstyle#1}$}}
\newcommand{\ldiagup}[1]%
       {\makebox[0cm]{${\scriptstyle#1}\uparrow\phantom{\scriptstyle#1}$}}
\newcommand{\rdiag}[1]%
       {\makebox[0cm]{$\phantom{\scriptstyle#1}\downarrow{\scriptstyle#1}$}}
\newcommand{\sediagr}[1]%
       {\makebox[0cm]{$\phantom{\scriptstyle#1}\searrow{\scriptstyle#1}$}}
\newcommand{\nediagr}[1]%
       {\makebox[0cm]{$\phantom{\scriptstyle#1}\nearrow{\scriptstyle#1}$}}
\newcommand{\rdiagup}[1]%
       {\makebox[0cm]{$\phantom{\scriptstyle#1}\uparrow{\scriptstyle#1}$}}
\newcommand{\swdiag}[1]%
       {\makebox[0cm]{$\phantom{\scriptstyle#1}\swarrow{\scriptstyle#1}$}}
\newcommand{\sediag}[1]%
       {\makebox[0cm]{${\scriptstyle#1}\searrow\phantom{\scriptstyle#1}$}}
\newcommand{\nediag}[1]%
       {\makebox[0cm]{${\scriptstyle#1}\nearrow\phantom{\scriptstyle#1}$}}

\newcommand{\doublearrowstack}[2]%
                      {{{{\scriptstyle#1}\atop{\textstyle\longrightarrow}}\atop{{\textstyle\longrightarrow}\atop{\scriptstyle#2}}}}
\newcommand{\rightleftarrowstack}[2]%
                      {{{{\scriptstyle#1}\atop{\textstyle\longrightarrow}}\atop{{\textstyle\longleftarrow}\atop{\scriptstyle#2}}}}
\newcommand{\leftrightarrowstack}[2]%
                      {{{{\scriptstyle#1}\atop{\textstyle\longleftarrow}}\atop{{\textstyle\longrightarrow}\atop{\scriptstyle#2}}}}

\newcommand{\overtoparrow}%
{\makebox[0cm]{\beginpicture
\setcoordinatesystem units <.8cm,.4cm> point at 0 0
\setplotarea x from -3 to 3, y from 0 to 1
\setquadratic
\plot -3 0 0 1 3 0 /
\put{\vector(3,-1){0}}[Bl] at 3 0
\endpicture}}

\newcommand{\underbottomarrow}%
{\makebox[0cm]{\beginpicture
\setcoordinatesystem units <.8cm,.4cm> point at 0 0
\setplotarea x from -3 to 3, y from 0 to 1
\setquadratic
\plot -3 1 0 0 3 1 /
\put{\vector(3,1){0}}[Bl] at 3 1
\endpicture}}

\newcommand{\ses}[5]%
{0\longrightarrow#1\stackrel{#2}{ \longrightarrow}#3\stackrel{#4}{
\longrightarrow}#5\longrightarrow0}

\newcommand{\dt}[6]%
{#1\stackrel{#2}{longrightarrow}#3 \stackrel{#4}{\longrightarrow}#5
\stackrel{#6}{\longrightarrow} #1[1]}  
 
\newcommand{\cat}[1]%
{(\mbox{\rm #1})}

\newcommand{\gm}{\Gamma }            
\newcommand{\lon }{\longrightarrow } 
\newcommand{\backl}{\mathbin{\vrule width1.5ex height.4pt\vrule height1.5ex}}
\newcommand{\per}{\backl}
\newcommand{\be}{\begin{eqnarray*}}
\newcommand{\ee}{\end{eqnarray*}}
\newcommand{\smalcirc}{\mbox{\tiny{$\circ $}}}

\newcommand{\half}{\frac{1}{2}}
\newcommand{\frakg}{\mathfrak g}

\newcommand{\ii}{i}

\newcommand{\Ind}{\mbox{Ind}}
\newcommand{\Ker}{\mbox{Ker}}

\title{Quantization of Pre-Quasi-Symplectic Groupoids and
their  Hamiltonian Spaces}
\author{Camille Laurent-Gengoux\\
Department of Mathematics\\
         Pennsylvania State University \\
         University Park, PA 16802, USA\\
{\sf email: laurent@math.psu.edu }
\\
AND\\
Ping Xu\thanks{Research partially supported by NSF grant  DMS03-06665. }\\
        Department of Mathematics\\
         Pennsylvania State University \\
         University Park, PA 16802, USA\\
{\sf email: ping@math.psu.edu }
}     
 
\date{}
\begin{document}
\sloppy
\maketitle

\centerline{{Dedicated to Alan Weinstein on the occasion of his 60th birthday}}

\begin{abstract}
We study the prequantization of pre-quasi-symplectic groupoids and their
Hamiltonian spaces using $S^1$-gerbes. We give a geometric description of the
integrality condition. As an application, we study the prequantization
of the quasi-Hamiltonian $G$-spaces of Alekseev--Malkin--Meinrenken.
\end{abstract}

\tableofcontents

\section{Introduction}


 Quasi-symplectic groupoids are natural generalizations of symplectic groupoids
 \cite{BCWZ, Xu}. The main motivation of \cite{Xu} in studying quasi-symplectic
 groupoids was to  introduce a single, unified momentum map theory in
 which ordinary Hamiltonian $G$-spaces, Lu's momentum maps of Poisson group
 actions, and the group-valued momentum maps of Alekseev--Malkin--Meinrenken
 can be understood under a uniform framework. An important feature of this
 unified theory is that it allows one to understand the diverse theories in
 such a way that techniques in one can be applied to the others.

It turns out that much of the theory of Hamiltonian spaces
of a symplectic groupoid can be generalized to 
 quasi-symplectic groupoids. In particular, one can perform
reduction and prove that $J^{-1}(\O)/\gm $ is a symplectic manifold, 
 where $\O\subset M$ is an orbit of the groupoid $\gm \toto M$.
More generally, one can introduce
the classical intertwiner space $\overline{X_{2}} \times_\gm X_1$
between two Hamiltonian $\gm$-spaces $X_1$ and $X_2$, generalizing the
notion studied by Guillemin-Sternberg \cite{GS} for
ordinary Hamiltonian $G$-spaces. One shows that
this is  a symplectic manifold (whenever it is a smooth manifold).
In particular, when $\gm$ is the AMM quasi-symplectic
groupoid \cite{BXZ, Xu}, this  reduced space describes the 
symplectic structure on the moduli space of flat connections
on a surface \cite{AMM}.

As is the case for symplectic groupoids \cite{Xu:90}, one
can introduce Morita
equivalence for quasi-symplectic groupoids. In particular, it has been
proven \cite{Xu} that (i) Morita equivalent quasi-symplectic groupoids give
rise to equivalent momentum map theories, in the sense that there
is a bijection between their Hamiltonian spaces; (ii) the
classical intertwiner space $\overline{X_{2}} \times_\gm X_1$
is independent of the Morita equivalence class of $\gm$. This
Morita invariance principle accounts for various well-known results
concerning the equivalence of momentum maps, including 
the Alekseev--Ginzburg--Weinstein linearization
theorem \cite{A, GW}
and the Alekseev--Malkin--Meinrenken equivalence
theorem   for group-valued momentum maps \cite{AMM} (see  \cite{Xu}
for details).

One important feature of Hamiltonian $G$-spaces is the
Guillemin--Sternberg theorem which states that ``$[Q, R]=0$'': quantization
commutes with reduction \cite{GS, Mein1}. One expects that
``$[Q, R]=0$'' should be a general guiding principle
for all momentum map theories. To carry out such a
quantization program, the first important step is
the construction of prequantum line bundles.
In this paper, we study the prequantization of Hamiltonian spaces for
quasi-symplectic groupoids. Our method uses the theory of
$S^1$-bundles and $S^1$-gerbes over a groupoid along with
their characteristic classes, as developed in \cite{BX, BX1}.
Roughly, our construction can be described as follows. 
A prequantization of  a  quasi-symplectic
groupoid $(\gm \toto M, \omega +\Omega )$ is  an $S^1$-central
extension $R\to \gm$  of the groupoid $\gm\toto M$ (or an $S^1$-gerbe
over the groupoid) equipped with  a pseudo-connection  having 
 $\omega +\Omega$ as  pseudo-curvature.
Such a prequantization exists if and only if
$\omega +\Omega$ is a de Rham integral $3$-cocycle
and $\Omega$ is exact (assuming that $\gm$ is
a proper groupoid). 
A prequantization
 of a Hamiltonian space is then  an $S^1$-bundle 
 $L$  over $R\toto M$ together with a compatible
pseudo-connection, where the $R$-action  on $L$ is $S^1$-equivariant.
 A prequantization of the symplectic intertwiner space
$\overline{X_{2}}\times_\gm X_1 $ can  be constructed
using these data.

Indeed one can show that
$R\backslash (\overline{L_{1}}\times_M  {L_2})$ is a prequantization
of the symplectic intertwiner space
 $ \overline{X_{2}}\times_\gm X_1 $, and the 
natural $1$-form on $ \overline{L_{1}}\times_M {L_2}$ induced
by the connection forms on $L_1$ and $L_2$ descends to
a prequantization connection on the quotient space 
$R\backslash (\overline{L_{1}}\times_M  {L_2})$.
When $\Omega $ is not exact, one must pass to a Morita equivalent
quasi-symplectic groupoid first. Then the Morita invariance principle
guarantees that the resulting quantization does not
depend on the  particular choice of Morita equivalent
quasi-symplectic groupoid.
As a special case, when $\gm$ is the AMM quasi-symplectic
groupoid, our construction yields the prequantization of
quasi-Hamiltonian $G$-spaces of
Alekseev--Malkin--Meinrenken and their symplectic reductions, 
and our  quantization condition 
coincides with that of Alekseev--Meinrenken \cite{AM}.

Quantization of Hamiltonian spaces for symplectic groupoids
was  studied in \cite{Xu:intert}.
Note that in the usual Hamiltonian case, since the symplectic
$2$-form defines a zero class in the third cohomology
group of the groupoid $T^{*}G\toto \frakg^*$, which
is the equivariant cohomology $H^3_G (\frakg^* )$,
gerbes do not  appear explicitly. However, for a general quasi-symplectic
groupoid (for instance  the AMM quasi-symplectic groupoid),
since the $3$-cocycle $\omega +\Omega$ may define a nontrivial class,
gerbes are inevitable in the construction.
Also note that no nondegeneracy condition is needed in
the quantization construction, so we drop this
assumption in the present paper to assure full generality.

This paper is organized as follows.
In Section 2, we review some basic results concerning
pre-quasi-symplectic groupoids and their Hamiltonian
spaces. In Section 3, we gather some important
results on $S^1$-bundles and $S^1$-central extensions.
We give  a simple formula for the index of an $S^1$-bundle
over a central extension in terms of the Chern class.
In Section 4, we introduce prequantizations  of pre-quasi-symplectic
groupoids and discuss compatible prequantizations
of their Hamiltonian spaces. Section  5 is devoted to the 
description of a geometric integrality condition
of pre-Hamiltonian $\gm$-spaces. The application to
quasi-Hamiltonian $G$-spaces is discussed.

Unless specified, by a groupoid in this paper, we always mean
a Lie groupoid whose orbit space is connected.
A remark  is in order concerning the terminology.
In \cite{BCWZ}, quasi-symplectic groupoids are
called presymplectic groupoids, where some ``non-degeneracy"
condition is assumed. Here we choose to use
the ``quasi" part of the terminology
to refer to the presence of a 3-form  and  to use ``pre-"
to men that ``non-degeneracy" is flexible.

Note that it would be interesting to investigate 
what notion of polarization would be relevant for the
next step of this quantization scheme.

Prequantization of symplectic groupoids was first studied
by Alan Weinstein and the second author in \cite{WX}, when
the second author was his PhD student. In the same paper,
$S^1$-central extensions  of Lie groupoids were also
 systematically investigated for the first  time.
Undoubtedly, Alan Weinstein's work and insights
have had a tremendous  impact on the development
of this  subject in the past two decades. It is
our great pleasure to dedicate this paper to him.

{\bf Acknowledgments.}
The second author would like to thank the Erwin Schr\"odinger Institute
and the University of Geneva for their hospitality
while work on this project was being done.
We would like to thank Anton Alekseev, Kai Behrend,
Eckhard Meinrenken,  and Jim Stasheff for useful discussions.

\section{Pre-Hamiltonian $\gm$-spaces and classical intertwiner spaces}

\subsection{Pre-quasi-symplectic groupoids and their pre-Hamiltonian spaces}

First, let us recall the definition of the 
 de Rham double complex of a 
Lie groupoid.
Let $\Gamma\toto M$ be a Lie groupoid.  Define for all
 $p\geq0$
$$\Gamma_{p}= 
\underbrace{\Gamma\times_{M}\ldots\times_{M}\Gamma}_{\text{$p$ times}}\,,$$
 {\em i.e.}, $\Gamma_{p}$ is the manifold of composable
sequences of $p$ arrows in the groupoid $\Gamma\toto M $ (and $\gm_0
=M$).
We have $p+1$ canonical maps $\gm_{p}\to \gm_{p-1}$ (each
leaving out
one of the $p+1$ objects involved in a sequence of composable arrows),
giving rise to a              diagram
\begin{equation}\label{sim.ma}
\xymatrix{
\ldots \gm_{2}
\ar[r]\ar@<1ex>[r]\ar@<-1ex>[r] & \gm_{1}\ar@<-.5ex>[r]\ar@<.5ex>[r]
&\gm_{0}\,.}
\end{equation}

Consider  the double complex $\Omega\upcom(\Gamma\lcom)$:
\begin{equation}
\label{eq:DeRham}
\xymatrix{
\cdots&\cdots&\cdots&\\
\Omega^1(\gm_{0})\ar[u]^d\ar[r]^\partial &
\Omega^1(\gm_{1})\ar[u]^d\ar[r]^\partial
&\Omega^1(\gm_{2})\ar[u]^d\ar[r]^\partial&\cdots\\
\Omega^0(\gm_{0})\ar[u]^d\ar[r]^\partial
&\Omega^0(\gm_{1})\ar[u]^d\ar[r]^\partial
&\Omega^0(\gm_{2})\ar[u]^d\ar[r]^\partial&\cdots
}
\end{equation}
Its boundary maps are $d:
\Omega^{k}( \gm_{p} ) \to \Omega^{k+1}( \gm_{p} )$, the usual exterior
derivative of differentiable forms and $\partial
:\Omega^{k}( \gm_{p} ) \to \Omega^{k}( \gm_{p+1} )$,  the
alternating sum of the pull-back maps of (\ref{sim.ma}).
We denote the total differential by $\delta = (-1)^pd+\del$.
The  cohomology groups of the total complex $C\upcom_{dR}(\Gamma\lcom)$
$$H_{dR}^k(\Gamma\lcom)= H^k\big(\Omega\upcom(\Gamma\lcom)\big)$$
are called the {\em de~Rham cohomology }groups of $\Gamma\toto M$.

\begin{defn}
A pre-quasi-symplectic groupoid is  a Lie groupoid  ${\gm}\toto{M}$
equipped with a $2$-form $\omega \in \Omega^2 (\gm )$ and
a $3$-form $\Omega\in \Omega^3 (M)$ such that
\begin{equation}
d \Omega =0, \ \   d \omega=\partial \Omega \ \ \mbox{and } \  \partial \omega
=0.
\end{equation}
In other words, $\omega+\Omega $ is  a 3-cocycle of the  total
de Rham complex of the groupoid  $\gm\toto M$.

A pre-quasi-symplectic groupoid $({\gm}\toto{M} , \omega+\Omega)$
is said to be {\em exact} if $\Omega$ is an exact $3$-form on $M$.
\end{defn}

A quasi-symplectic groupoid 
is a   pre-quasi-symplectic groupoid $({\gm}\toto{M} , \omega+\Omega)$
where $\omega$ satisfies certain non-degenerate condition \cite{BCWZ, Xu}.
Quasi-symplectic groupoids are natural generalization
of symplectic groupoids, whose momentum map theory
unifies   various  momentum
map theories, including the ordinary  Hamiltonian $G$-spaces,
Lu's momentum maps of Poisson group actions, and
group valued momentum maps of Alekseev--Malkin--Meinrenken.

\begin{defn}
\label{def:gm-space}
Given a pre-quasi-symplectic  groupoid  $({\gm}\toto{M}, \omega +\Omega )$,
a pre-Hamiltonian $\gm$-space is a (left) $\gm$-space $J: X\to M$
({\em i.e.}, $\gm$ acts on  $X$ from the left) 
with a compatible $2$-form $\omega_X \in \Omega^2 (X)$ such that:
\begin{enumerate}
\item $d \omega_X =J^* \Omega $;
\item the graph of the action $\Lambda=\{(r, x, rx)|t(r) =J(x)\}
 \subset \gm \times X \times \overline{X}$
(where $\overline{X}$ is the manifold $X$ endowed with the form
$-\omega_X$)
is isotropic with respect to the $2$-form 
$(\omega,  \omega_X, -{\omega_X  })$.
\end{enumerate}
\end{defn}
To illustrate the intrinsic meaning of the above
  compatibility condition, let us elaborate it in terms of groupoids.
Let $\gm \times_M X \toto X$ be the transformation groupoid
corresponding to the $\gm$-action,
 and, by abuse of notation,  $J: \gm \times_M X \to \gm$ the natural projection.
It is simple to see that
\begin{equation}
\label{eq:trans}
\xymatrix{
\gm \times_M X \ar@<-.5ex>[d]\ar@<.5ex>[d]\ar[r]^J
&
\gm \ar@<-.5ex>[d]\ar@<.5ex>[d]\\
X\ar[r]^J & M}
\end{equation}
  is a Lie groupoid homomorphism. Therefore
it induces a map, {\em i.e.},  the pull-back map, on the level of
the  de Rham complex:
$$J^* :\ \ \  \Omega^{\com}(\gm\lcom )\to \Omega^{\com}((\gm \times_M X)\lcom ). $$

\begin{prop} \cite{Xu}
\label{pro:2}
Let $({\gm}\toto{M}, \omega +\Omega )$ be a pre-quasi-symplectic 
 groupoid and $J: X\to M$
  a left $\gm$-space. A $2$-form $\omega_X\in \Omega^2 (X)$ is
compatible with the action
 if and only if
\begin{equation}
\label{eq:pr}
J^* ( \omega +\Omega )=\delta \omega_X.
\end{equation}
\end{prop}

\subsection{Classical intertwiner spaces}

Consider  a  pre-quasi-symplectic groupoid $({\gm}\toto{M}, \omega +\Omega )$,
and  pre-Hamiltonian $\gm$-spaces $(X_1~\stackrel{J_1}{\to}~M,~\omega_1~)$,
and $(X_2\stackrel{J_2}{\to} M, \omega_2)$.
Assume that $\gm \backslash (\overline{{X}_{2}}\times_M X_{1})$
is a smooth manifold, and denote by 
$$p: \overline{{X}_{2}}\times_M X_{1}\to 
\gm \backslash (\overline{{X}_{2}}\times_M X_{1})$$
the natural projection.
Note that $i^*(-{\omega_2}, \omega_1)$, where
$i: \overline{{X}_{2}}\times_M X_{1} \to X_{2} \times X_1$
is the natural embedding, is a closed $2$-form
on $\overline{{X}_{2}}\times_M X_{1}$.

\begin{prop}
\label{pro:intert}
The $2$-form $i^*(-{\omega_2}, \omega_1)$ descends to
a closed 2-form on $\gm \backslash (\overline{{X}_{2}}\times_M X_{1})$.
Therefore $\gm \backslash (\overline{{X}_{2}}\times_M X_{1})$ is
a presymplectic manifold.
\end{prop}               

To prove this proposition, we need a technical lemma.

\begin{lem}
\label{lem:red}
\cite{LTX}
Let $\gm\toto M$ be a Lie groupoid and $X\to M$ a 
left $\gm$-space. Assume that  $\gm \backslash X$
is a smooth manifold.  A 
differential form $\omega\in \Omega^* (X)$
descends to  a differential form on the quotient $\gm \backslash X$
if and only if $\partial \omega =0$, where
$\partial$ is with respect to the transformation
groupoid $\gm \times_M X \toto X$.
\end{lem}

\begin{pf1}
Note that the manifold ${{X}_{2}}\times_M X_{1}$ with the
momentum map $J: {{X}_{2}}\times_M X_{1} \to M$,
$J(x_2, x_1)=J_1 (x_1 )=J_2 (x_2 )$, is
 naturally a  $\gm$-space, where $\gm\toto 
M$ acts on ${{X}_{2}}\times_M X_{1}$  diagonally.
Then 
\be
\partial [ i^*(-{\omega_2}, \omega_1)]&=&i^* ( -\partial {\omega_2} ,\partial  \omega_1)\\
&=&i^* (-{J^*_2 \omega}  , J^*_1 \omega )\\
&=&\big((J_2\times J_1 )\smalcirc i)^* (-{\omega }  , \omega \big)\\
&=&0,
\ee
where $J_1, J_2$ and $i$ are respectively the
groupoid homomorphisms:

\begin{equation}
\label{eq:6}
\xymatrix{
\gm \times_M X_k \ar@<-.5ex>[d]\ar@<.5ex>[d]\ar[r]^{J_k}
&
\gm \ar@<-.5ex>[d]\ar@<.5ex>[d]\\
X_k\ar[r]^{J_k} & M} 
k=1, 2.
\end{equation}
and
\begin{equation}
\label{eq:7}
\xymatrix{
\gm \times_M (X_1  \times_M X_2)\ar@<-.5ex>[d]\ar@<.5ex>[d]\ar[r]^{i}
&
(\gm \times_M X_1) \times (\gm \times_M X_2) \ar@<-.5ex>[d]\ar@<.5ex>[d]\\
X_1 \times_M X_2\ar[r]^{i} & X_1 \times X_2} ,
\end{equation}
and $\partial [ i^*(-{\omega_2}, \omega_1)] $ and
$\partial \omega_k$, $k=1, \  2$,
are  with respect to the groupoids on the left-hand side
of Eqs.   (\ref{eq:7}) and (\ref{eq:6}), respectively.

The conclusion thus follows from Lemma \ref{lem:red}.
\end{pf1}

The presymplectic manifold
$\gm \backslash (\overline{{X}_{2}}\times_M X_{1})$ is called 
the {\em classical intertwiner space}, and is denoted
by $\overline{{X}_{2}}\times_\gm X_{1}$
for simplicity.
In particular, if   $({\gm}\toto{M}, \omega +\Omega )$  is a quasi-symplectic
groupoid, and  $(X_1\stackrel{J_1}{\to} M, \omega_1 )$
and $(X_2 \stackrel{J_2}{\to} M, \omega_2 )$ are Hamiltonian $\gm$-spaces, and
if  $J_1 :X_1  \to M$ and $J_2 : X_2\to M$ are clean,
then $\overline{{X}_{2}}\times_\gm X_{1}$ 
becomes a symplectic manifold. See \cite{Xu}
for details.

\section{$S^1$-bundles and $S^1$-central extensions}

In this section we   recall
 some basic results concerning $S^1$-bundles
 and $S^1$-central extensions over a groupoid. 
For details, consult \cite{BX, BX1, TXL}.

\subsection{Integral de Rham cocycles}

Let us recall some basic facts concerning singular homology.
For any manifold $N$, we denote by $(C_\com(N,{\mathbb Z}),d)$
the piecewise   smooth singular chain complex, 
and $Z_k(N, \zz )$ the space of smooth   $k$-cycles.
 For a smooth map $\phi: M \to N$,
we denote by $\phi_*$  both the chain map from $(C_\com(M,{\mathbb Z}),d)$
to $(C_\com(N,{\mathbb Z}),d)$ and the 
morphism of singular homology $ H_*(M,{\mathbb Z})\to H_*(N,{\mathbb Z}) $
induced by $\phi$.

For any Lie  groupoid $\Gamma \toto \gm_0$, 
consider  the double complex $C_\com(\Gamma_\com,{\mathbb Z})$:
$$ 
\begin{array}{cccccc}
\cdots&  &\cdots&  &\cdots \\
\downarrow d & & \downarrow d & & \downarrow d\\
C_1(\gm_0,{\mathbb Z}) & \stackrel{\partial}{\leftarrow} & C_1(\gm_1,{\mathbb Z})& \stackrel{\partial}{\leftarrow} & C_1(\gm_2,{\mathbb Z}) \\
\downarrow d & & \downarrow d & & \downarrow d\\
C_0(\gm_0,{\mathbb Z}) & \stackrel{\partial}{\leftarrow} & C_0(\gm_1,{\mathbb Z})& \stackrel{\partial}{\leftarrow} & C_0(\gm_2,{\mathbb Z}), \\
\end{array}$$
 where $ \gm_0 =M$, and $\partial:C_k(\Gamma_p,{\mathbb Z})
 \to C_k(\Gamma_{p-1},{\mathbb Z}) $
is the alternating sum of the chain  maps induced by the face maps.
We denote the total differential by $\delta=(-1)^p d + \partial $.
Its homology will be denoted by $H_k(\Gamma_\com , {\mathbb Z}) $.
By $Z_k(\Gamma_\com , {\mathbb Z}) $ we denote
the space of  $k$-cycles
and by $[C] \in H_* (\Gamma_\com , {\mathbb Z}) $ the class of a given cycle
$C$. Note that $C_k(\Gamma_p, \zz)$ is the free Abelian group generated by the
piecewise smooth maps $\Delta_k\to\Gamma_p$.

The construction above can be carried out in  exactly the
same way   replacing ${\mathbb Z}$ by $ {\mathbb R} $.
The corresponding  homology groups are denoted by 
 $H_k(\Gamma_\com , {\mathbb R}) $.
According to the universal-coefficient formula
(see, for example, \cite{Spanier}),
there is a canonical isomorphism 
$$H_k(\Gamma_\com , {\mathbb R})\simeq
H_k(\Gamma_\com , {\mathbb Z})\otimes_{\mathbb Z} {\mathbb R} .$$ 

There is a natural pairing between  $C_\com(\Gamma_\com,{\mathbb R})$
and $C_{dR}^\com(\Gamma_\com)$ given as follows.
For any generator $C: \Delta_k \to \Gamma_p$ in
 $C_\com(\Gamma_\com,{\mathbb Z})$,

\begin{equation} \label{eq:pairing}
 \langle C,\omega\rangle = \left\{\begin{array}{cc} \int_{\Delta_k} C^* \omega &  \hspace{0.2cm} 
\mbox{if $\omega \in \Omega^k(\Gamma_p) $}\\
 0 &\hspace{0.2cm}  \mbox{otherwise.} \\
 \end{array}\right.
\end{equation}
For simplicity, we will denote this pairing by $\int_C \omega$.
With this notation, the pairing satisfies the following identities:
\be \label{eq:relations}&&  \int_{d C} \omega = \int_C d\omega \\ 
   &&    \int_{\partial C} \omega = \int_C \partial\omega   \\ &&
  \int_{\delta C} \omega = \int_C \delta\omega    \\ &&     \ee

Moreover, if $\phi: G \to H$ is a groupoid homomorphism,
then for any $C \in C_\com(G_\com,{\mathbb R})$
and  $\omega \in C^\com_{dR}(H_\com)$ 
    \begin{equation} \label{eq:eq:morphisme} \int_{\phi_*(C)}\omega= \int_C
 \phi^* \omega  .  \end{equation}
%

The following result is standard (see, for example, Proposition 6.1
in  \cite{Dupont}).

\begin{prop}
 \label{prop:degenerescence}
The pairing $  H_k(\Gamma_\com,{\mathbb R}) \otimes
H^k_{dR}(\Gamma_\com )\to {\mathbb R}, \ \
([C],[\omega]) \to \int_C \omega $  is non-degenerate.
\end{prop}
                              
Let
 
\begin{equation}\label{def:hk} 
Z^k_{dR}(\Gamma_\com,{\mathbb Z}) =
\{\omega \in   Z^k_{dR}(\Gamma_\com) | \  \int_C \omega \in {\mathbb Z} \mbox{ for any cycle $C \in Z_k(\Gamma_\com,{\mathbb Z})$ }    \}  . \end{equation}
Elements in $Z^k_{dR}(\Gamma_\com,{\mathbb Z})$ are called
{\em integral de Rham cocycles}, or simply integral cocycles. 


\subsection{$S^1$-bundles and $S^1$-central extensions}

In this subsection, we recall some basic notations and results
concerning $S^1$-bundles and $S^1$-central extensions over a Lie 
groupoid. For details, see \cite{BX, BX1}.

\begin{defn}
Let $\Gamma\toto M$ be a Lie groupoid. A (right) {\em $S^1$-bundle }over
$\Gamma\toto M$ is a (right) $S^1$-bundle $P$ over $M$, together with
a (left) action of $\Gamma$ on $P$ which respects the $S^1$-action
({\em i.e.}, we have $(\gamma \cdot x) \cdot t=\gamma \cdot (x \cdot t )$,
for all $t\in S^1$ and  all compatible pairs $(\gamma,x)\in
\Gamma\times_{M} P$).
\end{defn}

Let $Q\toto P$ denote the corresponding  transformation groupoid 
$\Gamma\times_{M} P \toto P$.
There is a natural groupoid homomorphism  $\pi$ from $Q\toto P$ to
 $\Gamma\toto M$.  Of course, $Q$ is an $S^1$-bundle over $\Gamma$.

A {\em pseudo-connection}
 is a $1$-cochain $\theta \in C^1_{dR}(Q\lcom)$,
where  $\theta\in \Omega^1(P)$ is  a connection 1-form for the $S^1$-bundle
$P\to M$.
 One checks that $\delta\theta \in C^2_{dR}(Q\lcom)$
descends to a 2-cocycle in $Z^2_{dR}(\Gamma\lcom )$.
 In other words, there exist
unique $\omega\in \Omega^1(\Gamma)$ and $\Omega\in\Omega^2(M)$ such
that 
$$\delta\theta =\pi\upst(\omega+\Omega). $$
Then $\omega+\Omega$ is called the {\em pseudo-curvature}, which 
is an integral $2$-cocycle.
Its class $[\omega+\Omega]\in H^2(\Gamma\lcom , \zz )$ is
called the {\em Chern class} of the  $S^1$-bundle $P$.

\begin{prop}
\label{pro:3.2}
\cite{BX, BX1}
 Let $\gm \toto M$ be  a proper Lie groupoid.
Assume that $\omega+\Omega\in  \Omega^1 (\gm )\oplus \Omega^2(M)\subset
C^2_{dR} (\Gamma\lcom)$ is an integral
2-cocycle.  Then there exists an $S^1$-bundle $P$
over $\Gamma\toto M$ and a pseudo-connection  $\theta 
\in \Omega^1(P) $ for
the bundle $P\to M$  whose pseudo-curvature
equals  $\omega+\Omega$.
\end{prop}

\begin{defn}
Let ${\gm}\toto {M}$ be a Lie groupoid. An {\em $S^1$-central
extension }of $\gm\toto M$ consists of

1) a  Lie groupoid ${R}\toto {M}$, together with a morphism of Lie
groupoids $(\pi,\id):[R\toto M]\to[\Gamma\toto M]$,

2) a left $S^1$-action on $R$, making $\pi:R\to \Gamma$ a (left)
principal $S^1$-bundle.
\noindent These two structures are compatible in the sense that
$(s\cdot x) (t\cdot y)=st \cdot (xy )$,
for all  $ s,t \in S^{1}$ and $(x, y) \in R\times_MR $.
\end{defn}

Given a central extension $R$ of $\Gamma\toto M$,  a {\em pseudo-connection} 
 is a $2$-cochain $\theta +B\in  C^2_{dR}(R\lcom)$, where
 $\theta\in \Omega^1(R)$  is a connection
1-form for the bundle $R\to \Gamma$ and $B\in \Omega^2 (M)$.
It is simple to check that $\delta (\theta +B)$ descends
to a 3-cocycle  in $Z^3 (\gm\lcom)$, {\em i.e.},
$$\delta (\theta +B)=\pi^* (\eta +\omega +\Omega)$$
 for some $\eta +\omega +\Omega\in Z^3 (\gm\lcom)$.
Then $\eta +\omega +\Omega$ is an integral cocycle
in $Z^3_{dR} (\gm\lcom, \zz )$, and is  called the {\em pseudo-curvature}.
Its class $[\eta +\omega +\Omega]\in H^3 (\gm\lcom , \zz )$
 is called the Dixmier-Douady class of $R$.

\begin{prop}
\cite{BX, BX1}
Assume that $\gm \toto M$ is
a proper Lie groupoid.
Given any  3-cocycle  $\eta +\omega+\Omega\in Z^3_{dR}(\Gamma\lcom)$
 such that

1) $[\eta+\omega+\Omega]$ is integral, and

2) $\Omega$ is exact,

\noindent there exists a groupoid $S^1$-central extension $R\toto M$ of the
groupoid $\gm \toto M$, and a pseudo-connection $\theta +B \in
\Omega^1 (R)\oplus \Omega^2 (M)$
such that its pseudo-curvature equals
$\eta+\omega+\Omega$.
\end{prop}

\subsection{Index of an $S^1$-bundle over a central extension}
Let $R \stackrel{\pi}{\to} \Gamma \toto M$ be an
$S^1$-central extension,
and $S^1 \to L \stackrel{p}{\to} M$  
a principal $S^1$-bundle over the groupoid
$R \toto M$ with Chern class $[L]\in H^1 (R\lcom, S^1 )$.
The  example below will be useful in the future.

\begin{numex}
 Consider, for any $k \in {\mathbb Z}$, the
principal $S^1$-bundle $B_k: S^1 \to \cdot$ over $S^1 \toto \cdot$,
where  the groupoid $S^1 \toto \cdot$ 
acts on $B_k$ by
 $$\lambda \cdot z= \lambda^k z \hspace{1cm} \forall \lambda \in S^1 \toto \cdot 
\hspace{0.4cm}\mbox{and} \hspace{0.4cm}\forall z \in S^1 \to \cdot. $$

It is well-known that  $ H^1(S^1_\com,S^1) \simeq {\mathbb Z}$.
Under this isomorphism,
the class $[B_k]$ is simply equal to  $k$.

It is also simple to see that the Chern class of  $B_k$
  can be represented by
\begin{equation}\label{eq.Chernclass}  k  \frac{dt}{2\pi} \in
 Z^{1} (S^1) \subset Z^2 ( (S^1)_\com ),
 \end{equation}
where $\frac{dt}{2\pi}$ is the normalized Haar measure on $S^1$.
\end{numex}

For any $m \in M$, there exists a groupoid homomorphism
 $f_m$ from $S^1 \toto \cdot$ to $R \toto M$
defined by
\begin{equation}\label{eq:plonge} f_m(\lambda) = \lambda \cdot 1_m
 \hspace{1cm} \forall \lambda \in S^1 ,
\end{equation}
where $1_m \in R$ is the unit element over $m \in M$.

This  homomorphism induces a map
\begin{equation} \label{eq:index} 
f^*_m \,: \,  H^1(R_\com,S^1) \to  H^1(S^1_\com,S^1) \simeq {\mathbb Z}.
\end{equation}

 For a  principal $S^1$-bundle $L$ over $R \toto M$,
we define its {\em index}
by 
$$ \Ind_m (L) = f^*_m([L]) \, \,  \in H^1(S^1_\com,S^1 ) \simeq {\mathbb Z}.$$

We  list some of its important  properties below.

\begin{prop}\label{prop:index}
Let $R \stackrel{\pi}{\to} \Gamma \toto M$ be an
$S^1$-central extension,
and $S^1 \to L \stackrel{p}{\to} M$
an (right) principal $S^1$-bundle over the groupoid
$R \toto M$.  Then
\begin{enumerate}
\item the index is characterized by the relation 
$$ f_m(\lambda)\cdot l =l \cdot  \lambda^{\Ind_m(L)},
   \, \,  \, \, \forall
\lambda \in S^1, \  l \in p^{-1}(m),$$
where the  dot  on the left hand side denotes
the $R$-action on $L$, while  the  dot  on
the right hand side refers to the $S^1$-action
on $L$;
\item for any $m \in M$, the pull-back $ f_m^* L$ is isomorphic
 to $B_{\Ind_m(L)}$;
 \item $\Ind_m(L)$ is constant on the groupoid orbits;
\item   $\Ind_m(L)$ is constant on any connected component of $M$; and
\item if $\gm \backslash M $ is path connected, then the index $\Ind_m(L)$ 
is independent of  $m \in M$.
\end{enumerate} 
\end{prop}
\begin{pf}
1) and 2)
Let $ l$ be any point  in the fiber  $ L_{m} =p^{-1} (m)$.
For any $\lambda \in S^1$, there exists a unique $\phi(\lambda) \in S^1 $
such that 

\begin{equation}\label{eq:caracterisation}
 f_m(\lambda)\cdot l =l \cdot  \phi(\lambda).
\end{equation}
The map $\lambda \to \phi(\lambda)$ does not depend on the choice of $l$
in the fiber $p^{-1}(m)$ and
is a group homomorphism from $S^1$ to  $S^1$. Therefore, it is 
  of the form
$\phi( \lambda)= \lambda^k$ for some $k \in {\mathbb Z}$.

Now $f_m^*([L]) \in H^1(S^1_\com,S^1)$ is the Chern class associated to
the pull-back of $L$ by $f_m$.
On the other hand, according to Eq.~(\ref{eq:caracterisation}),   
$f_m^* L$
is isomorphic (as a principal $S^1$-bundle over $S^1 \toto \cdot$) to $B_k$.
Therefore $k =\Ind_m(L)$. and Eq.~(\ref{eq:caracterisation}) implies 
\begin{equation}\label{eq:caracterisation2}
 f_m(\lambda)\cdot l =l \cdot  \lambda^{\Ind_m(L)},
   \, \,  \, \, \forall l \in p^{-1}(m).
\end{equation}
This proves 1) and 2).

3) For  any $\gamma \in R$ with $s(\gamma)=n$ and $t(\gamma)=m$,
we have   $\gamma 1_m=1_{ n} \gamma$.
  It follows from Eq.~(\ref{eq:plonge}) that 
$ \gamma f_m(\lambda)   =  f_{ n} (\lambda) \gamma $.
Now for any $ l \in p^{-1}(m)$,   we have
$(\gamma f_m(\lambda)) \cdot l= ( f_{ n} (\lambda) \gamma  )\cdot l $.
On the one hand, we have 
\begin{equation}\label{eq:caract}
 (\gamma f_m(\lambda)) \cdot l =(\gamma \cdot l )   \cdot \lambda^{\Ind_m(L)},
\end{equation}
and
\begin{equation}\label{eq:caract2}
 ( f_{ n} (\lambda) \gamma  )\cdot l=
f_{ n} (\lambda) ( \gamma  \cdot l)=
(\gamma \cdot l) \cdot \lambda^{\Ind_{n}(L)}.
\end{equation}
From Eqs.~(\ref{eq:caract}) and (\ref{eq:caract2}), it
follows that $\Ind_m(L)=\Ind_{n}(L)$.

4) It is clear from Eqs.~(\ref{eq:caracterisation2})
that  $\Ind_m(L)$ depends continuously on $m \in M$.
Since it is a ${\mathbb Z}$-valued function, we have $\Ind_m(L)=\Ind_n(L) $
for any pair of points  $(m,n) \in M\times M$ that are 
in the same connected component of $M$.

5) follows from 3) and 4) immediately. 
\end{pf}

\subsection{Index and Chern class}

From now on, we  will assume that the space  of orbits
$M / \Gamma$ is  path-connected,  and denote  the
 index of $L$  simply by $\Ind(L)$.
Therefore we have a group homomorphism:
$$\Ind(L) :  H^1(R_\com,S^1 ) \to \zz.$$

From the commutativity of the diagram
$$ \begin{array}{ccc} H^1(R_\com,S^1 ) &
\stackrel{\ii}{\to}& H^2(R_\com,{\mathbb Z} )\\
\downarrow  & & \downarrow \\
 H^1(S^1_\com,S^1) & \stackrel{\ii}{\to} & H^2(S^1_\com,{\mathbb Z}) \\
 \end{array},$$
we see that 
$\Ind(L) $ factors through $H^2(R_\com,{\mathbb Z} )\to \zz$.
In the following proposition, we give an explicit formula
for $\Ind(L) $ in terms of  the Chern class.

\begin{prop}
\label{prop:integrale}
 Assume that $L\to M$ is a principal
$S^1$-bundle over $R \toto M$ with the Chern
class $[\theta +\omega ] \in H^2_{dR} (R\lcom , \zz )$,
 where $R \stackrel{\pi}{\to} \Gamma \toto M$
is an $S^1$-central extension,
 and $\theta + \omega \in \Omega^1(R) \oplus \Omega^2(M)$.
 Then the index of $L$ is given by
   $$  \Ind(L)= \int_{\pi^{-1}(\epsilon(m))} \theta  ,$$ 
where $\epsilon : M \to \Gamma$ is the unit map.
\end{prop}
\begin{pf} 
Let $L'$ be the pull-back of the principal $S^1$-bundle
$L$ via  the  homomorphism  $f_m:S^1 \to R$.
The Chern class of $L' $ is the pull-back of the Chern class
of $L$, {\em i.e.}, the class defined by
$f_m^* \theta + f_m^* \omega \in C^2_{dR} (S^1_\com) $.
Since $f_m^* \omega$ is a $2$-form over a point, it vanishes
and therefore  the Chern class of $L' $  is
 represented by $f_m^* \theta \in \Omega^1(S^1)$.

By Proposition \ref{prop:index}, $L'$ is  isomorphic to $B_{\Ind(L)}$. 
According to Eq.~(\ref{eq.Chernclass}), 
 the identity $f_m^*~\theta~=~\Ind(L)~\frac{dt}{2\pi}~+~\delta g
= ~\Ind(L)~\frac{dt}{2\pi}~+d g $ holds
for some function $g \in C^{\infty}( S^1,{\mathbb R})$.

Now since $f_m$ is a bijection from $S^1$ to $\pi^{-1}(\epsilon(m)) $,
we have 
  $$  \int_{\pi^{-1}(\epsilon(m))} \theta=  \int_{S^1}f_m^* \theta.$$

 Therefore
   $$  \int_{\pi^{-1}(\epsilon(m))} \theta=  \int_{S^1}f_m^* \theta=\Ind(L) \int_{S^1} \frac{dt}{2\pi}+ \int_{S^1} d g=\Ind(L)  .$$
\end{pf}

Recall that a line bundle $L\to M$
over $R\toto M$ is called
a $(\gm, R)$-twisted line bundle if $\ker \pi \cong M\times S^{1}$
 acts on $L$  by scalar multiplication, where $S^1$ is identified with
the unit circle in $\cc$ \cite{TXL}.
The following corollary is an immediate consequence of 
Proposition \ref{prop:integrale} and  Proposition \ref{prop:index}.

\begin{cor}\label{cor:twisted} 
Under the hypotheses of Proposition \ref{prop:integrale},
$L\to M$ defines a $(\gm, R)$-twisted line bundle if and only if
$$\int_{\pi^{-1}(\epsilon(m))} \theta =1.$$
\end{cor}



\section{Prequantization of classical intertwiner spaces}

\subsection{Compatible prequantizations}

\begin{defn}
A prequantization of a  pre-quasi-symplectic  groupoid
$({\gm}\toto{M}, \omega +\Omega )$ consists 
of an $S^1$-central extension $R\stackrel{\pi}{\to} \gm \toto M$
together with a pseudo-connection $\theta + B\in 
\Omega^1 (R) \oplus \Omega^2 (M)$ such that
\begin{equation}
\label{eq:5}
\delta(\theta+B)= \pi\upst(\omega+\Omega). 
\end{equation}
\end{defn}

According to Proposition \ref{pro:3.2}, if ${\gm}\toto{M}$ is
a proper Lie groupoid, a prequantization exists if and
only if $({\gm}\toto{M}, \omega +\Omega )$ is
exact and $\omega+\Omega$ is an integral  3-cocycle.
A pre-quasi-symplectic  groupoid
$({\gm}\toto{M}, \omega +\Omega )$  is said to be {\em integral}
if $\omega+\Omega$ is an integral  cocycle.

\begin{defn}
\label{def:compatible}
Let $({\gm}\toto{M}, \omega +\Omega )$ be an
exact pre-quasi-symplectic  groupoid,
 and $(R\to \gm \toto M, \theta + B )$ a prequantization.
Assume that  $(X\stackrel{J}{\to} M, \omega_X)$ is
 a pre-Hamiltonian $\gm$-space. A compatible prequantization
of $X$   consists of  an $S^1$-bundle $\phi: L\to X$ with a connection
$1$-form $\theta_L \in \Omega^1 (L) $ such that
\begin{enumerate}
\item  $\tilde{J}=J\smalcirc  \phi   : L\lon M$ is
a left $R$-space  and the action satisfies:
$$ (s\cdot \kappa ) (t\cdot x)=st \cdot (\kappa x ), $$
for all  $ s,t \in S^{1}$ and $(\kappa, x) \in R\times_MX $
a compatible pair;
\item the $1$-form $ (\theta ,\theta_{L} , -\theta_{L} )\in
\Omega^{1}(R\times L\times \overline{L})$ vanishes on  the graph of the action
$$\Xi =\{(\kappa, l, \kappa l)|\kappa\in R, l\in L \ \mbox{compatible
pairs}\};$$ and
\item $d\theta_L =\phi^* (J^*B-\omega_X)$.
\end{enumerate}
\end{defn}

Note that the  second condition above is equivalent to  saying 
that
$(R\times L\times \overline{L})/T^{2}\stackrel{p}{\lon }\gm \times X \times X $
with $p([ \kappa, l, m ])=(\pi (\kappa ),\phi (l),\phi (m))$ is a
flat $S^1$-bundle with the   connection $\bar{\Theta }$, which 
is the $1$-form  on   $ (R\times L\times \overline{L})/T^{2}$  naturally
induced from $\Theta =(\theta ,\theta_{L} , -\theta_{L} )\in
 \Omega^{1}(R\times L\times L)$ (see \cite{Xu:intert}).

\begin{numex}\label{ex:t*G}
If $\gm$ is the symplectic groupoid
$(T^*G \toto \frakg^* , \omega )$, where
$\omega\in \Omega^2 (T^*G )$ is the
canonical cotangent symplectic 2-form, a prequantization
of $\gm$ can be taken to be $R\cong T^*G \times S^1\to T^*G$,
the trivial $S^1$-bundle and $\theta =\theta_{T^*G}+ dt$,
where $\theta_{T^*G}\in \Omega^1 (T^*G )$ is the Liouville
$1$-form and $t$ is the natural coordinate on $S^1$.
A Hamiltonian $\gm$-space is a Hamiltonian $G$-space  $J: X\to 
\frakg^*$ in the usual sense. It is simple to see
that a compatible pre-quantization is a $G$-equivariant
prequantization of $X$, which always exists when
$G$ is connected and simply connected \cite{GS}.
\end{numex}

More generally, the following result was proved in
\cite{Xu:intert} (the theorem was stated for  the symplectic
case, but it is valid for  the presymplectic case as well).

\begin{prop}
Let $(\gm \toto M, \omega )$ be an   $s$-connected and $s$-simply
connected pre-symplectic  groupoid, and $(X\stackrel{J}{\to} M, \omega_X )$
a pre-Hamiltonian  space. If both $\omega $  and $\omega_X$ represent
 integral cohomology classes in $H^2_{dR}(\Gamma)$ and
$H^2_{dR}(X)$ respectively, then there exists a 
compatible prequantization.
\end{prop}

For a given pre-quasi-symplectic  groupoid
$({\gm}\toto{M}, \omega +\Omega )$ and
a prequantization $(R\to \gm \toto M, \theta+ B )$,
let $\gm\times_M X\toto X$ be the transformation
groupoid as in Eq. (\ref{eq:trans}). By pulling  back the 
central extension $R\to \gm \toto M$ via $J$,  one obtains a
central extension of groupoids $R\times_M X\to \gm\times_M X\toto X$.
Here $R\times_M X$ is again a transformation groupoid,
 where $R$ acts on $X$ by projecting $R$ to $\gm$ and
using the given $\gm$-action on $X$.

By abuse of notation, we still use $J$ to denote the
projection $R\times_M X\to R$. Therefore we have the following
homomorphism of $S^1$-central extensions of groupoids:
\begin{equation}
\xymatrix{
R \times_M X \ar[d]\ar[r]^J
&
R \ar[d]\\
\gm \times_M X \ar@<-.5ex>[d]\ar@<.5ex>[d]\ar[r]^J
&
\gm \ar@<-.5ex>[d]\ar@<.5ex>[d]\\
X\ar[r]^J & M}
\end{equation}

\begin{numrmk}
Note that Proposition \ref{pro:2} implies that the  Dixmier-Douady class of 
$R\times_M X\to \gm \times_M X\toto X$ vanishes.
If $\gm\toto M$ is a proper groupoid, so is $\gm \times_M X\toto X$.
Therefore $R\times_M X\to \gm \times_M X\toto X$ defines a 
trivial gerbe. According to Proposition 4.2  of \cite{BX},
there exists an $S^1$-bundle $E\to X$ such that
$R\times_M X \cong s^* E\otimes t^*\overline{E}$
as a  central extension.
\end{numrmk}

\begin{prop}
\label{pro:5.4}
Let $({\gm}\toto{M}, \omega +\Omega )$ be an
exact pre-quasi-symplectic  groupoid,
 and $(R\to \gm \toto M, \theta  +B )$ its prequantization.
Assume that  $(X\stackrel{J}{\to} M, \omega_X)$ is
 a pre-Hamiltonian $\gm$-space.
Then $(L\stackrel{\phi}{\to} X, \theta_L)$ is a     compatible prequantization
of $X$  if and only if (the associated line bundle of)
$\phi: L\to X$ is a twisted line bundle over
$R\times_M X\to \gm \times_M X\toto X$ with   the  pseudo-connection
and the pseudo-curvature  being given by $\theta_L$ and 
$J^*\theta+(J^*B-\omega_X) \in \Omega^1(R\times_M X)\oplus \Omega^2 (X)$,
respectively.
\end{prop}
\begin{pf}
Given a compatible prequantization $L\stackrel{\phi}{\to} X$,
define an action of $R\times_M X\toto X$ on $L$ by
$(\kappa, \phi (l))\cdot l= \kappa l$, where $\kappa\in R$ and
$l\in L$ are compatible pairs. It is simple to check that
all the compatibility conditions are satisfied so that
$L\stackrel{\phi}{\to} X$ is a twisted line bundle over 
the central extension $R\times_M X\to \gm \times_M X\toto X$. It is
simple to see that the corresponding transformation 
groupoid $(R\times_MX)\times_X L \toto L$ is
isomorphic to the transformation groupoid 
$ R\times_M L\toto L$.
Moreover,  it is simple to see that Condition (2) of 
Definition \ref{def:compatible} implies that
\begin{equation}
\label{eq:19}
\partial \theta_L = \phi^* J^* \theta,
\end{equation}
where, by abuse of notation, we use $\phi$ to denote the
Lie groupoid homomorphism:
\begin{equation}
\label{eq:Q1}
\xymatrix{
R\times_M L \ar@<-.5ex>[d]\ar@<.5ex>[d]\ar[r]^\phi &
R\times_MX \ar@<-.5ex>[d]\ar@<.5ex>[d]\\
L\ar[r]^\phi & X}
\end{equation}
and $\partial \theta_L $ is with respect to the
groupoid $ R\times_M L\toto L$.
Therefore we have
$$\delta \theta_L=\partial \theta_L +d\theta_L
= \phi^* (J^*\theta +J^*B-\omega_X). $$

The converse can be proved by working backwards. 
\end{pf}

As an immediate consequence, we have

\begin{cor}
\label{cor:int}
Under the hypotheses of Proposition \ref{pro:5.4}
and assuming  that $\gm\toto M$ is proper,  for a pre-Hamiltonian $\gm$-space
$(X\stackrel{J}{\to} M, \omega_X)$,
a compatible prequantization exists if and only if
$ J^*(\theta +B)-\omega_X $ is an integral 2-cocycle
in $Z^2_{dR}((R\times_M X)\lcom, \zz)$.
\end{cor}
\begin{pf}
One direction is obvious by Proposition \ref{pro:5.4}.

For the other direction, note that Proposition
\ref{pro:2}  implies that $J^*(\theta+B)-\omega_X$ is
always a $2$-cocycle since
$$\delta (J^*(\theta+B)-\omega_X )=
J^* \delta (\theta+B)-\pi^* \delta \omega_X=
J^* \pi^* (\omega+\Omega)-\pi^*J^* (\omega+\Omega)=0.$$
Here we have used Eqs. (\ref{eq:5}) and (\ref{eq:6}).
If $ J^*(\theta+B)-\omega_X $ is an integral cocycle in
$Z^2_{dR}((R\times_M X)\lcom, \zz)$, according to Proposition \ref{pro:3.2},
 there exists an $S^1$-bundle $L\to X$
over $R\toto X$ and a pseudo-connection  $\theta_L \in \Omega^1(L) $
 whose pseudo-curvature equals to $J^*\theta+(J^*B-\omega_X)$.
According to Corollary \ref{cor:twisted}, one sees that  (the associated line bundle of) $L$ is indeed
a twisted line bundle over $R\times_M X\to \gm \times_M X\toto X$.
Then $L\to X$ is a compatible prequantization by Proposition \ref{pro:5.4}. 
\end{pf}

\subsection{Prequantization of classical intertwiner spaces}

We are now ready to state the main theorem of this section.

\begin{them}
\label{thm:main}
Let $({\gm}\toto{M}, \omega +\Omega )$ be an
exact pre-quasi-symplectic  groupoid,
 and $(R\to \gm \toto M, \theta+ B )$ a prequantization.
Assume that  $(X_k\stackrel{J_k}{\to} M, \omega_k), \ k=1, 2$,  are
pre-Hamiltonian $\gm$-spaces, $\gm \toto M$
acts freely on $\overline{X_2}\times_M X_1$, and
 $\overline{X_2}\times_\gm X_1=\gm \backslash(\overline{X_2}\times_M X_1)$
is a smooth manifold.
Let
$(L_k\stackrel{\phi_k}{\to} X_k, \theta_k)$, 
be a compatible prequantization of $X_k$ for $k=1,2$.
Then
$$\phi:
R \backslash ( L_2 \times_M \overline{L_1})\to \overline{X_2}\times_\gm X_1, \ 
 \phi [l_2, l_1]=[\phi_2 (l_2 ), \phi_1 (l_1)],  $$
 with the $S^1$-action $\lambda \cdot [l_2, l_1]=[\lambda \cdot l_2,   l_1],
\ \lambda \in S^1,$
 is an $S^1$-principal bundle. Moreover,
 $i^* (\theta_2, - \theta_1)$
descends to a connection $1$-form on 
$R \backslash ( {L_2}\times_M \overline{L_1})$, which
defines a prequantization of the
classical intertwiner space $\overline{X_2}\times_\gm X_1$.
Here $i: L_2\times_M \overline{L_1}\to L_2\times \overline{L_1}$ is the natural embedding.
\end{them}
\begin{pf}
One checks directly that $\phi:
R \backslash ( {L_2}\times_M L_1)\to \overline{X_2}\times_\gm X_1$
is an $S^1$-bundle. Now let $R\toto M$ act on ${L_2}\times_M L_1$
diagonally.
We have
$$\partial i^* (\theta_2, -\theta_1)
=  i^* (\partial \theta_2, -\partial \theta_1 )
=i^* (\phi_2^*J_2^* \theta ,- \phi_1^* J_1^* \theta )=0.$$
Hence $i^* (\theta_2, - \theta_1)$ descends to a 1-form
on the quotient space $R \backslash ( \overline{L_2}\times_M L_1)$,
which can be easily seen to be a connection 1-form.
Now 
$$d(i^* (\theta_2, - \theta_1)) =
i^*(d\theta_2, - d\theta_1)=
i^* ( \phi_2^*(J_2^*B-\omega_2), \phi_1^* (J_1^*B-\omega_1) )
=i^* (\phi_2 \times \phi_1 )^* (-\omega_2, \omega_1 ),$$
where in the last equality  we used
the relation   $J_1\smalcirc \phi_1=J_2\smalcirc \phi_2$
on ${L_2}\times_M L_1$.
Here $\phi_i$ and $J_k$, $k=1, 2$ are groupoid homomorphisms:
\begin{equation}
\label{eq:transf}
\xymatrix{
R \times_M L_k \ar@<-.5ex>[d]\ar@<.5ex>[d]\ar[r]^{\phi_k}
& R\times_M X_k\ar@<-.5ex>[d]\ar@<.5ex>[d]\ar[r]^{J_k}
&R \ar@<-.5ex>[d]\ar@<.5ex>[d]\\
L_k\ar[r]^{\phi_k} & X_k \ar[r]^{J_k} & M}
\end{equation}
and  $i$ is the  groupoid homomorphism:
\begin{equation}
\label{eq:grouhomo}
\xymatrix{
R \times_M (L_2  \times_M L_1)\ar@<-.5ex>[d]\ar@<.5ex>[d]\ar[r]^{i}
&
(R \times_M L_2) \times (R \times_M L_1) \ar@<-.5ex>[d]\ar@<.5ex>[d]\\
L_2 \times_M L_1\ar[r]^{i} & L_2 \times L_1}.
\end{equation}

This  completes the proof.
\end{pf}


\subsection{Morita equivalence}\label{sec:morita}

\begin{defn}
Pre-quasi-symplectic groupoids $({G}\toto G_0 , \omega_G +\Omega_G )$
and $({H}\toto H_0 , \omega_H +\Omega_H )$ are  said to
be {\em  Morita equivalent}
if there exists  a  Morita equivalence bimodule
 $G_{0}\stackrel{\rho }{\leftarrow} X
\stackrel{\sigma}{\rightarrow} H_{0}$ between the Lie groupoids $G$ and $H$,
together with a $2$-form $\omega_X\in \Omega^2 (X)$ such that
$(X\stackrel{\rho\times \sigma}{\to}G_0 
\times \overline{H_0}, \omega_X )$ is  a 
pre-Hamiltonian $G\times \overline{H}$-space, 
where the $G\times \overline{H}$-action  on $X$ is
given by
$(g, h)\cdot x=gxh^{-1}$ for all compatible triples 
$ g\in G, \ h\in H$ and $  x\in X$.
\end{defn}               

One easily checks that this is indeed an equivalence
relation among pre-quasi-symplectic groupoids.

Let $Q\toto X$ be the transformation groupoid
$$Q:(G\times \overline{H})\times_{(G_0\times  \overline{H_0} )} X\toto X. $$
Then the  natural projections $\pr_1: Q\to G$ and
$\pr_2: Q\to  \overline{H}$ are groupoid homomorphisms.
As an immediate consequence of Proposition \ref{pro:2}, 
 we have the following   identity 
$$\pr_1^* (  \omega_G +\Omega_G )- \pr_2^*  (  \omega_H +\Omega_H )
=\delta \omega_X . $$
Note that the axioms of Morita equivalence of
Lie groupoids assure that,  as groupoids, 
$Q\cong G[X]$ and $Q\cong H[X]$ (see the proof of Proposition 4.5 
of \cite{Xu}), where $G[X]\toto X$ and $H[X]\toto X$ are the pull-back
groupoids of $G$ and $H$ using $\rho$ and $\sigma$, respectively.

Recall that for a given Lie groupoid $\gm \toto M$,
two cohomologous 3-cocycles $\omega_i +\Omega_i\in
\Omega^2 (\gm )\oplus \Omega^3 (M )$, $i=1, 2$,
are said to differ by a  {\em gauge transformation of
the first type} if 
$$(\omega_1+\Omega_1)-(\omega_2+\Omega_2) =\delta B$$
for some $B\in \Omega^2 (M)$. 

By a {\em Morita morphism} from  the  pre-quasi-symplectic
groupoid  $(\gm'\toto M', \omega'+\Omega' )$ 
to $(\gm\toto M, \omega+\Omega )$, we mean  a Morita  morphism of
the Lie groupoid  $p: \gm'\to \gm$ ({\em i.e.} $\gm'$ is isomorphic to
the pullback groupoid $\gm [M' ]\toto M'$) such that
$\omega'+\Omega' $ and $p^*\omega +p^* \Omega $ differ
by a gauge transformation of the first type.

The following result gives a more intuitive explanation of Morita
equivalence.

\begin{prop}
Two pre-quasi-symplectic groupoids are Morita equivalent if
 and only if there exists a third pre-quasi-symplectic
groupoid together with a Morita morphism to each of them.                           \end{prop} 

\begin{cor} 
For two Morita equivalent pre-quasi-symplectic groupoids,
if one is integral, so is the other.
\end{cor} 

Therefore Morita equivalence   induces
an equivalence relation among integral pre-quasi-symplectic groupoids.

One of the most important features of Morita equivalent
pre-quasi-symplectic groupoids is the following

\begin{them}
\label{thm:morita}
Suppose that $({G}\toto G_0 , \omega_G +\Omega_G )$
and $({H}\toto H_0 , \omega_H +\Omega_H )$ are Morita equivalent
pre-quasi-symplectic groupoids with
  an equivalence bimodule $G_{0}\stackrel{\rho }{\leftarrow} X\stackrel{\sigma}{\rightarrow} H_{0}$.  Then:
\begin{enumerate}
\item Corresponding to
any pre-Hamiltonian $G$-space $J_F:F \to G_0$, there is a unique
(up to isomorphism)
pre-Hamiltonian $H$-space $J_E:E \to H_0$
 such that $F$ and $E$ are a pair of related pre-Hamiltonian spaces
 and vice versa.
\item Let $J_{F_i} :F_i \to G_0$, $i=1, 2$, be pre-Hamiltonian $G$-spaces
and $J_{E_i}: E_i\to H_0$, $i=1, 2$, their related  pre-Hamiltonian
$H$-spaces. Then $\overline{F_2}\times_G F_1$ and
$\overline{E_2} \times_H E_1$ are diffeomorphic as presymplectic manifolds
(in the sense that if  one is smooth so is the other).
\end{enumerate}                                       
\end{them} 
\begin{pf}
This was proved in \cite{Xu} for quasi-symplectic groupoids 
and their Hamiltonian spaces. One can prove this theorem
in a similar fashion (in fact in a simpler way by using
Proposition \ref{pro:intert}). We will leave the details
to the reader.
\end{pf} 

We now can introduce   Morita equivalence for the prequantization
of pre-quasi-symplectic groupoids.

\begin{defn}
Let $({G}\toto G_0 , \omega_G +\Omega_G )$
and $({H}\toto H_0 , \omega_H +\Omega_H )$ be Morita equivalent 
integral exact pre-quasi-symplectic groupoids with
an  equivalence bimodule
 $(G_{0}\stackrel{\rho }{\leftarrow} X
\stackrel{\sigma}{\rightarrow} H_{0}, \ \omega_X)$.
We say their  prequantizations   $(R_{G}\to {G}\toto G_0 , \theta_G+B_G)$
and $(R_{H}\to {H}\toto H_0 , \theta_H+B_H)$  are Morita
equivalent if   $X$  admits a compatible prequantization
$(Z\to X, \theta_Z)$ with respect to
 the prequantization of the  pre-quasi-symplectic groupoid
 $(R_G\times \overline{R_H})/S^1\to G\times \overline{H}\toto
G_0\times \overline{H_0}, (\theta_G+\overline{\theta_H})+(B_G+\overline{B_H}))$.
\end{defn}

It is simple to see that $G_{0}{\leftarrow} Z
{\rightarrow} H_{0}$ is an equivalence bimodule of
central extensions in the sense of Definition 2.11 \cite{TXL}.

\begin{numrmk}
\begin{enumerate}
\item Note that  prequantizations can be Morita equivalent as
central extensions, but not Morita equivalent 
as prequantizations. The former one simply means
that they correspond to isomorphic $S^1$-gerbes,
and up to a torsion, are determined by their  Dixmier-Douady classes.
\item It would be interesting to investigate the following question:
given two Morita equivalent pre-quasi-symplectic
groupoids and a prequantization of one of them,
is it possible to construct a Morita equivalent 
 prequantization for the other pre-quasi-symplectic
groupoid? 
\end{enumerate}  
\end{numrmk}

A useful  feature of Morita equivalence is that it gives 
a recipe which allows us to   construct compatible
prequantizations.

\begin{them}
For Morita equivalent  prequantizations of pre-quasi-symplectic
groupoids, there is an equivalence of categories of
compatible prequantizations of pre-Hamiltonian spaces.
\end{them}
\begin{pf}
Let $({G}\toto G_0 , \omega_G +\Omega_G )$
and $({H}\toto H_0 , \omega_H +\Omega_H )$ be Morita equivalent
integral exact pre-quasi-symplectic groupoids with
an  equivalence bimodule
 $(G_{0}\stackrel{\rho }{\leftarrow} X
\stackrel{\sigma}{\rightarrow} H_{0}, \ \omega_X)$,  and 
$(R_{G}\to {G}\toto G_0 , \theta_G+B_G)$
and $(R_{H}\to {H}\toto H_0 , \theta_H+B_H)$  be Morita
equivalent  prequantizations  given by
$(Z\to X, \theta_Z)$.  
Assume that $J:F \to G_0$ is a pre-Hamiltonian $G$-space
and $(L\to F, \theta_L)$ a compatible  prequantization.
It is known that the corresponding  pre-Hamiltonian
$H$-space is  $E:=\overline{X}\times_G F\stackrel{J'}{\to} H_0$,
where $J' : E\to H_{0}$ and the  $H$-action on $E$ are defined by
$ J' ([x,f])=\sigma (x) $ and $ h\cdot [x, f]=[x\cdot h^{-1}, f]$, respectively.

Let $L'=Z\times_{R_G}\overline{L}$. Then it is clear that $L'$ is
an $S^1$-bundle   over $E$, and $R_H$ acts on $L'$ equivariantly.
It is simple to check that $i^* (\theta_Z, -\theta_L)$,
where $i: Z\times_{G_0}L\to Z\times L$, descends to a 
$1$-form on the quotient space $Z\times_{R_G}\overline{L}$
which is  indeed  a connection $1$-form $\theta_{L'}$ on $L'$.
It is routine to check that $(L'\to E, \theta_{L'})$ is
a compatible prequantization of the pre-Hamiltonian
$H$-space $J' : E\to H_{0}$.

The inverse functor can be constructed in a similar fashion. 
\end{pf}

\begin{numrmk}
The above theorem indicates a useful method which enables one
to transform  
  prequantizations
of  Hamiltonian $LG$-spaces to prequantizations    of quasi-Hamiltonian
$G$-spaces of AMM and vice-versa.  The latter is understood as a compatible
prequantization   corresponding to  the quasi-symplectic
groupoid $(G\times G)[{\mathcal U}]\toto \coprod U_i$, which
is the pull-back quasi-symplectic groupoid of the AMM quasi-symplectic groupoid
using an open covering  ${\mathcal U}= (U_i)_{i\in I}$
 of $G$ (see \cite{M},  for instance, 
for an explicit construction). It is known that
 $(G\times G)[{\mathcal U}]\toto \coprod U_i$  is
Morita equivalent to the  symplectic groupoid $(LG\times L\frakg \toto { L\frakg},
\omega_{LG \times L\frakg} )$ according to Proposition 4.26 \cite{Xu}.
The  question is therefore boiled down to the construction of
a compatible prequantization of the  Morita equivalence Hamiltonian bimodule.
\end{numrmk}

\section{Integral pre-Hamiltonian $\gm$-spaces}
The main purpose of this section is  to give a geometric
integrality condition which guarantees the
existence of  a prequantization of a 
pre-Hamiltonian $\gm $-space.
 
\subsection{Integrality condition}

\begin{lem}\label{lem:condition} Let $J: G \to H$ be a 
groupoid homomorphism.
By $\Ker(J_*)$ we denote
 the kernel of $J_*: H_2(G_\com,{\mathbb Z})
\to H_2(H_\com,{\mathbb Z})$.
Let $\omega \in Z^2_{dR}(G_\com )$.
 The following conditions are equivalent:
\begin{enumerate}
\item there exists   $\Xi \in Z^2_{dR}(H_\com)$ such that
 $$ \omega + J^* \Xi  \in  Z^2_{dR}(G_\com, {\mathbb Z});$$
\item for any $C \in Z_2(G_\com, {\mathbb Z})$ with $ [C] \in \Ker(J_*)$, 
we have
 $$\int_C \omega \in \zz   .$$
\end{enumerate}
\end{lem}
\begin{pf}
 1) $\Rightarrow$ 2).
By   definition,
we have  for any $C \in Z_2(G_\com,{\mathbb Z})$, 
\begin{equation}\label{eq:moduloZ1} \int_C ( \omega + J^* \Xi) \in \zz.
    \end{equation}
From Eq.~(\ref{eq:relations}), we  also have
  $$ \int_C ( \omega + J^* \Xi)= \int_C \omega + \int_{J_*(C)} \Xi .  $$
If $[J_*(C)]=0$, {\em i.e.},  $J_*(C)=\delta D$ for some $D \in 
C_3(H_\com,{\mathbb Z})$, then
$$ \int_C ( \omega + J^* \Xi)= \int_C \omega + \int_{\delta D} \Xi 
=\int_C \omega + \int_{D} \delta \Xi= \int_C \omega $$
since $\delta \Xi=0$.
Therefore
$ \int_C \omega  \in \zz $.

2) $\Rightarrow$ 1). Since there exists  
a $\zz$-submodule 
 $\hH$ in $H_2(G_\com,{\mathbb Z})$ 
such that $H_2(G_\com,{\mathbb Z})= \hH \oplus \Ker(J_*) $,
the ${\mathbb Z}$-map 
$$f:\Ker(J_*) \to {\mathbb Z}, \ f([C])=\int_C \omega, \ \ \forall
 [C]\in \Ker(J_*) ,$$ 
can be extended to a ${\mathbb Z}$-map $\tilde{f}: H_2(G_\com,{\mathbb Z}) 
\to \zz$.  According to  Proposition \ref{prop:degenerescence},
there exists  $\omega' \in  Z^2_{dR}(G_\com)$
such that 
 $$\tilde{f}([C])=\int_C \omega' ,\ \ \ \forall [C] \in H_2(G_\com,{\mathbb Z}).$$
By Eq. (\ref{def:hk}), $\omega' $ is an  integral cocycle
in $Z^2_{dR}(G_\com,{\mathbb Z})$.
Moreover,  we have
\begin{equation} \label{eq:ortho} \int_C (\omega' -\omega)=0,  \ \ \ \ 
\forall C \in Z_2 (G_\com , {\mathbb Z})  \ \mbox{ such that }
[C]\in \Ker(J_*) .\end{equation}
Since $J_*: H_2 (G_\com , {\mathbb R}  ) \to
 H_2 (H_\com , {\mathbb R}  )$ is dual to
 $J^*: H^2_{dR} (H_\com)
\to H^2_{dR} (G_\com  )$, we have
$\Ker(J_*) ^{\perp}=\mbox{Im} (J^*)$.
Therefore $[\omega' -\omega] = J^* [\Xi] $ for some $\Xi \in Z^2_{dR}(\gm_\com)$.
 This proves 1).
\end{pf}

\begin{defn}
Let $(\Gamma \toto M, \omega + \Omega)$ be a
 pre-quasi-symplectic groupoid.
A pre-Hamiltonian $\gm$-space $ (X \to M, \omega_X) $
is said to satisfy the {\em integrality condition}  if
for any $C \in Z_2((\Gamma \times_M X)_\com,{\mathbb Z})$ and  any
$D \in C_3(\Gamma_\com,{\mathbb Z}) $ 
  \begin{equation}
\label{eq:integrable}
 \delta D = J_*(C) \hspace{1cm} \Rightarrow \hspace{1cm}  \int_C \omega_X - \int_{D} (\omega + \Omega) \in {\mathbb Z}.
\end{equation} 
In this case, we also say that the pair $(\omega_X , \omega + \Omega)$
satisfies the integrality condition.
\end{defn}

\begin{numrmks}
\begin{enumerate}
\item By taking $C=0$, Eq. (\ref{eq:integrable}) implies that
$\int_{D} (\omega + \Omega) \in {\mathbb Z}, \ \ \forall D\in Z_3 (\Gamma_\com,{\mathbb Z}) $. That is, $\omega + \Omega$ must be an integral 3-cocycle
and therefore 
$(\Gamma \toto M, \omega + \Omega)$ must be
 an integral  pre-quasi-symplectic groupoid.

\item If $\omega + \Omega$ is  a 3-coboundary $\delta K$,
then the integrality condition is equivalent to
\begin{equation}     
\int_C (\omega_X-J^*K)\in \zz, \ \ \ \forall C \in 
Z_2((\Gamma \times_M X)_\com,{\mathbb Z})  \ \mbox{such that } J_*[C]=0.
\end{equation} 
\end{enumerate}
\end{numrmks}   

From now on, we shall always assume that
 $(\Gamma \toto M, \omega + \Omega)$ is an integral
 pre-quasi-symplectic groupoid. The following lemma
indicates that it is sufficient to require that
both sides of  Eq. (\ref{eq:integrable})
hold for a single representative $(C, D)$ in every  class of $\Ker (J_*)$.

\begin{lem}
\label{lem:simpl}
Let $(\Gamma \toto M, \omega + \Omega)$ be an integral
 pre-quasi-symplectic groupoid.
A pre-Hamiltonian $\gm$-space $ (X \to M, \omega_X) $
 satisfies the  integrality condition  if   
and only if for any
  class  $\Cc\in \Ker J_*$, there exists 
  $ C \in Z_2((\Gamma \times_M X)_\com,{\mathbb Z})$ and
  $D \in  C_3(\Gamma_\com,{\mathbb Z})$
 with $\Cc =[C]$ and  $J_*(C)=\delta D$
 such that 
  $$ \int_C \omega_X - \int_{D} (\omega + \Omega) \in {\mathbb Z}.  $$
\end{lem}
\begin{pf}
 Let   $C' \in Z_2((\Gamma \times_M X)_\com,{\mathbb Z})$ and
$D' \in C_3(\gm_\com,{\mathbb Z}) $ be any pair
satisfying $J_*(C')=\delta D'$.  Then $[C'] \in \Ker (J_*)$.
By  assumption,
there exists a pair $(C,D)$ such that $[C]=[C']$,
 $J_*(C)=\delta D$, and $\int_{C} \omega_X -\int_{D}(\omega+\Omega )
\in \zz$.
Assume that
 $C=C'+ \delta E $ for some $E \in C_3((\Gamma \times_M X)_\com,{\mathbb Z})$.
Then we have
\be
&&\int_{C'} \omega_X - \int_{D'}(\omega + \Omega) -
\big(\int_{C} \omega_X - \int_{D}(\omega + \Omega) \big)\\
&=&-\int_{\delta E} \omega_X 
  + \int_{D} (\omega+ \Omega)-
  \int_{D'}(\omega + \Omega) \\
&=&-\int_{E} \delta \omega_X + \int_{D} (\omega+ \Omega)-
  \int_{D'}(\omega + \Omega) \\ 
&=&-\int_{E} J^*(\omega + \Omega) + \int_{D} (\omega+ \Omega)-
  \int_{D'}(\omega + \Omega) \\   
&=&\int_{D-J_*(E)-D'  } (\omega + \Omega)  .
\ee
Since $\delta(D-J_*(E)-D' )=J_*(C-C'-\delta E) =0$ and
$\omega + \Omega$ is an integral cocycle, 
it follows that $\int_{D-J_*(E)-D'  } (\omega + \Omega)  \in \zz$.
This completes the proof. 
\end{pf}

Now assume that $(\Gamma \toto M,\omega+\Omega)$
is an integral  exact  pre-quasi-symplectic groupoid, and  
$R \to \Gamma \toto M$  is a prequantization.
Let $\theta + B \in  \Omega^1 (R) \oplus \Omega^2 (M)$ 
be  a pseudo-connection satisfying   Eq.~(\ref{eq:5}).
 In order to fix the notation, recall
 that we have  the following commutative diagram of
groupoid homomorphisms
\begin{equation}\label{eq:not}\begin{array}{ccc}
  R \times_M X   & \stackrel{J}{\to}& R  \\
   \downarrow \pi &  & \downarrow \pi \\
\gm \times_M X   & \stackrel{J}{\to}& \gm 
\end{array}\end{equation}
where the horizontal arrows are  projections.

\begin{lem}
\label{lem:7.5}
Assume that  $C' \in Z_2((R \times_M X)_\com,{\mathbb Z})$
satisfies $J_*( {C'})= k Z + \delta D'$ for some $D' \in C_3(R_\com, {\mathbb Z})$.
Let $C=\pi_*(C')$  and  $D= \pi_*(D')$. Then
\begin{equation}
\label{eq:machin4}   
\int_C \omega_X -\int_D (\omega+ \Omega)=
k+\int_{C'} \big(\omega_X - J^*(\theta + B)\big),
  \end{equation}           
where $Z \in Z_1 (R\lcom ,\zz)$ is the 1-cycle defined by   Eq.~(\ref{defZm}).
\end{lem}
\begin{pf}
First, since $\pi: R \times_M X \to \Gamma \times_M X$ reduces
to the identity map when being restricted to the unit spaces,
we have
 
   \begin{equation}\label{eq:machin1} \int_{C'} \omega_X = \int_{\pi_*({C'})} \omega_X  = \int_{C} \omega_X.\end{equation}

Now by Eq.~(\ref{eq:relations}), we have
   \begin{equation}
\label{eq:machin2}
  \int_{C'} J^* (\theta + B) = \int_{J_*({C'})} (\theta + B)
=k \int_{Z} (\theta + B) + \int_{\delta D'}(\theta +B ) .
 \end{equation}

According to  Lemma \ref{lem:precisions}, 
$ \int_{Z} (\theta + B)= \int_{Z} \theta  =1 $.
 Therefore
 \be 
  \int_{C'} J^* (\theta + B) &=&k+  \int_{\delta D'}(\theta +B ) \\
&=&  k+\int_{D'}\delta(\theta +B ) \hspace{1.5cm} \mbox{(by Eq.~(\ref{eq:5}))}  \\
&=&  k+ \int_{D'}\pi^*(\omega + \Omega)   \hspace{1.5cm}
\mbox{(by Eq.~(\ref{eq:relations}))}   \\
&=&  k+ \int_{D}(\omega + \Omega). 
\ee
Hence it follows that
$$   \int_C \omega_X -\int_D (\omega+ \Omega)=
k+\int_{C'} (\omega_X - J^*(\theta + \Omega)).$$
\end{pf}

The following proposition gives  a useful characterization of
integrality condition.

\begin{prop}
\label{thm:int} 
Let $(\Gamma \toto M, \omega + \Omega)$ be an integral
 pre-quasi-symplectic groupoid, and $(R\to \gm \toto M, \theta+B)$
a prequantization.  
Assume that $ (X \stackrel{J}{\to} M, \omega_X) $
is a pre-Hamiltonian $\gm$-space.
Then the following conditions are equivalent.
\begin{enumerate}
\item There exists a 2-cocycle $\Xi \in Z^2_{dR}(\Gamma_\com)$ such that
$$ \omega_X - J^* (\theta +B) -J^* \pi^* \Xi \in Z^2_{dR}\big((R \times_M X)_\com,{\mathbb Z}\big) .$$
\item For any cycle $C' \in Z_2( (R \times_M X)_\com,{\mathbb Z})$
such that $[C']\in \Ker(\pi_* \smalcirc J_*)$, we have
  $$ \int_{C'} \big( \omega_X - J^* (\theta +B)  \big) \in \zz .$$
\item The pair $(\omega_X,\omega+\Omega)$ satisfies the integrality condition.
\end{enumerate}
\end{prop}
\begin{pf}
 1) $\iff$ 2).  follows from  Lemma \ref{lem:condition}.

 2) $\Rightarrow$ 3).
 Any class in $ \Ker J_* \subset H^2((\gm  \times_M X)_\com,{\mathbb Z})$
can be represented by a 2-cocycle of the form
  $C=\pi_*(C')$ where $C' \in Z_2((R \times_M X)_\com,{\mathbb Z})$.
Then $[C']$ is in the kernel of $J_* \smalcirc \pi_* = \pi_* \smalcirc J_*$.
It thus follows that $[J_*(C')] \in \Ker(\pi_*) $.
By Lemma \ref{lem:precisions}, $[J_*(C')]= k [Z]$ for some 
 $k \in {\mathbb Z}$ where $Z \in C_1(R_\com,{\mathbb Z})$
is defined by Eq.~(\ref{defZm}). In other words, there exists
$D' \in C_3(R_\com, {\mathbb Z})$ such that
    $J_*( {C'})= k Z + \delta D' $.
Let  $D= \pi_*(D')$. One can  easily see 
that $\delta D=\pi_*(J_*(C'))=J_*(C)$.
Then by Lemma \ref{lem:7.5}, we have
$\int_C \omega_X -\int_D (\omega+ \Omega) \in \zz$.
By Lemma \ref{lem:simpl}, this implies that the pair $(\omega_X,\omega+\Omega)$
 satisfies the integrality condition.
   
   3) $\Rightarrow$ 2).
Let $C' \in Z_2((R \times X)_\com ,{\mathbb Z})$ be any  cycle whose class is
in the kernel of $ \pi_* \smalcirc J_*$.
Since $[J_*(C')] \in \Ker(\pi_*)$,   
 Lemma \ref{lem:precisions} implies that
there exists $k \in {\mathbb Z}$ and $D' \in C_3(R_\com, \zz )$
such that 
  \begin{equation}\label{eq:machin6}   J_*(C')= k Z + \delta D'. \end{equation}
Therefore,  by Eq. (\ref{eq:machin4}), we have
$  \int_{C'}(\omega_X-J^* (\theta+B))= -k+\int_{C}\omega_X -
 \int_{D} (\omega + \Omega)$,
where $C=\pi_*(C')$ and $D=\pi_* (D')$.
By  applying $\pi_*$ to  Eq.~(\ref{eq:machin6}), one finds
that  $J_*(C)= \delta D$.  Since $(\omega_X, \omega+\Omega)$
satisfies the integrality condition,  it thus follows that
$  \int_{C'}(\omega_X-J^* (\theta+B))\in  \zz$.
\end{pf}

As an immediate consequence, we  obtain the following main
result of the section.

\begin{them}\label{cor:quant}
Let $(\Gamma \toto M,\omega+\Omega)$
be  an   exact proper  pre-quasi-symplectic groupoid,
and  $(X \stackrel{J}{\to} M, \omega_X)$  a
pre-Hamiltonian $\gm$-space. Then there exists a compatible 
prequantization  $R\to \Gamma \toto M$ and $L\to X$
if and only if the pair $(\omega_X, \omega+\Omega)$   satisfies the
integrality  condition of Eq. (\ref{eq:integrable}).
\end{them} 
\begin{pf}
Assume that $(R\to \Gamma \toto M, \theta +B)$ and $(L\to X, \theta_L)$
are a pair of compatible prequantizations. By 
Corollary \ref{cor:int}, we have
 $\omega_X - J^* (\theta +B)  \in Z^2_{dR}((R \times_M X)_\com,{\mathbb Z})$.
Hence $(\omega_X, \omega+\Omega)$   satisfies the
integrality  condition according to  Proposition \ref{thm:int}.

Conversely, assume that $(\omega_X, \omega+\Omega)$   satisfies the
integrality  condition.
Then $\omega+\Omega$ must be an integral cocycle.
 Let $(R\to \Gamma \toto M, \theta +B)$
be a prequantization, which always exists since 
$\gm$ is proper. Again according to Proposition \ref{thm:int},
there exists  a 2-cocycle
 $\Xi \in Z^2_{dR}(\Gamma_\com)$ such that
$ \omega_X - J^* (\theta +B) -J^* \pi^* \Xi \in 
Z^2_{dR}((R \times_M X)_\com,{\mathbb Z}) $. Since 
$\gm$ is proper, $\Xi$ is cohomologous to $\alpha +B_0$, where
$\alpha \in \Omega^1 (\gm)$ is a closed 1-form and  $ B_0\in \Omega^2 (M)$.
Then $\theta' +B':
=(\theta +\pi^*\alpha) +(B+B_0)$ is clearly  also a
pseudo-connection and $ \omega_X - J^* (\theta' +B')  
 \in Z^2_{dR}((R \times_M X)_\com,{\mathbb Z}) $. 
From  Corollary \ref{cor:int}, it follows that
$(X \stackrel{J}{\to} M, \omega_X)$
admits  a  compatible prequantization
$(L\to X, \theta_L)$.
\end{pf}

\subsection{Integral quasi-Hamiltonian $G$-spaces}

In this subsection, $G$ is a  connected and
 simply-connected  compact  Lie group
and $1$ denotes the unit of $G$.
We intend  to study the case where $\gm $
is the AMM quasi-symplectic groupoid.

Assume that $X$ is  a $G$-space. 
There is  a natural map $i:H_2(X,{\mathbb Z}) \to H_2\big((G \times X)_\com,{\mathbb Z}\big) $
induced by the inclusion  $C_2(X,{\mathbb Z}) \subset
 C_2((G \times X)_\com,{\mathbb Z})$.
The following lemma indicates that $i$ is in fact  an isomorphism.

\begin{lem}\label{lem:trans-cycle} 
If $G$ is a connected and simply-connected Lie group, then the map
  $$i: H_2(X,{\mathbb Z}) \to H_2((G \times X)_\com,{\mathbb Z})$$
is an isomorphism.
\end{lem}
\begin{pf}
This is a standard result. For completeness,
we sketch a  proof below.
Let $G \to EG \to BG$ be the usual $G$-bundle
over the classifying space $BG$ 
and $X_G = G \backslash (EG \times X) $.
We have the fibration $ G \to EG\times X \to X_G$.

The second term of the homology Leray-Serre spectral sequence
 is $E^2_{p,q}=H_p(X_G, {\mathcal H}_q(G, \zz) )$,
{\em i.e.,} the homology of $X_G$ with local coefficients in 
${\mathcal H}_q(G, \zz)  $
(see \cite{McCleary}).
Since $G$ is simply-connected, we have $H_1(G,{\mathbb Z})=H_2(G,{\mathbb Z})=0$
and $E^2_{p,q}$ has the following form for $0 \leq p \leq 3$
and $ 0 \leq q \leq 2$: 
\begin{equation} \label{eq:spectre} \begin{array}{ccccc}
 \vdots & \vdots & \vdots & \vdots &  \\
  0  & 0 & 0 & 0 & \cdots \\
0  & 0 & 0 & 0   & \cdots \\
H_0(X_G,{\mathbb Z})  & H_1(X_G,{\mathbb Z}) & H_2(X_G,{\mathbb Z}) &
H_3(X_G,{\mathbb Z})& \cdots  \\
\end{array}
\end{equation}

 According to Leray-Serre
theorem, this spectral sequence converges 
to $H_*(EG \times X,{\mathbb Z})$.
It is clear from Eq.~(\ref{eq:spectre})
that, in particular, we have
 $$ H_2(EG \times X,{\mathbb Z})\simeq H_2(X_G,{\mathbb Z}).$$
Since $EG$ is contractible,  we get
$$  H_2(X_G,{\mathbb Z})\simeq H_2(X,{\mathbb Z}).$$
The lemma now follows  from the well-known 
isomorphism   $H_2((G \times X)_\com,{\mathbb Z})\simeq H_2(X_G,{\mathbb Z}) 
$.
\end{pf}

 Since $H_2(\Gg^*,{\mathbb Z})=0$,  Lemma \ref{lem:trans-cycle} implies
that
\begin{equation} \label{ref:sing-hom-grou-1} H_2\big((T^* G)_\com,{\mathbb Z}\big)=0 .\end{equation}
Since any simply-connected Lie group  $G$ satisfies
$H_2(G,{\mathbb Z})=0$, we also  have 
\begin{equation} \label{ref:sing-hom-grou-2}
 H_2\big((G \times G)_\com ,{\mathbb Z}\big)=0 .\end{equation}

Recall that the  AMM quasi-symplectic groupoid
is $(G \times G \toto G, \omega+\Omega )$ \cite{BXZ, Xu},
where $G$  is  a compact  Lie group equipped with an 
ad-invariant non-degenerate symmetric
bilinear form $(\cdot , \cdot )$.
Here $G \times G \toto G$ is the transformation groupoid,
where $G$ acts on itself by conjugation, and $\omega$ and 
$\Omega $ are defined as follows.

Following  \cite{AMM},  we denote by $\theta $ and $\bar{\theta}$
the left and right Maurer-Cartan forms on $G$ respectively, i.e.,
$\theta =g^{-1}dg$ and $\bar{\theta}=(dg ) g^{-1}$.
Let $\Omega \in \Omega^{3}(G)$ denote the bi-invariant
3-form on $G$ corresponding to the Lie algebra 3-cocycle
$\frac{1}{12}(\cdot , [\cdot , \cdot ])\in \wedge^3 \frakg^*$  {\em i.e.}
\begin{equation}
\label{eq:chi}
\Omega =\frac{1}{12}(\theta , [\theta , \theta ])
=\frac{1}{12}(\bar{\theta}, [\bar{\theta} , \bar{\theta} ]),
\end{equation}
and  $\omega \in \Omega^2 (G\times G )$  the $2$-form
\begin{equation}
\label{eq:quasi}
\omega|_{(g, x)} =-\half [(Ad_{x} \pr_1^* \theta , \pr_1^* \theta )
+(\pr_1^* \theta , \pr_2^{*}(\theta +\bar{\theta} ))],
\end{equation}
where $(g, x)$ denotes the coordinate in $G\times G$, and
$\pr_1$ and $\pr_2: G\times G\to G$ are natural projections.
It is known that  $\omega +\Omega$ is a integral 3-cocycle.

A  triple $(X,\omega_X,J)$, where $X$ is a manifold,
 $\omega_X$ is a $G$-invariant $2$-form on $X$ and  $J:X\to G$
is  a smooth map,
is a quasi-Hamiltonian $G$-space in the sense of \cite{AMM} if
 \begin{enumerate} 
\item[(B1)]\label{B1} the differential  of $\omega_X$ is given by: 
\begin{equation*}  \label{Gclosed} 
d\omega_X=  J^*\Omega;
 \end{equation*}                     
\item[(B2)]\label{B2}
the map $J$ satisfies
\begin{equation*}
\hat{\xi} \per \omega_X =\half J^*  (\xi, \theta +\bar{\theta });
\end{equation*} and
\item[(B3)]\label{B3}
at each $x\in X$, the kernel of $\omega_X$ is given by
\begin{equation*} 
\ker \omega_X= \{ \hat{\xi}(x)|\ \xi\in\ker(\Ad_{J(x)}+1) \},
 \end{equation*}
\end{enumerate}
where $\hat{\xi}$ is  the vector field on $X$ associated to
the infinitesimal action of $\xi \in \Gg$ on $X$.

It is known \cite{Xu} that
these conditions are equivalent to  $(X\stackrel{J}{\to}G,  \omega_X)$
being a   Hamiltonian  $\gm$-space,  where $\gm$ is the
 AMM quasi-symplectic groupoid $(G \times G \toto G,\omega+\Omega)$.
In this case, the integrality can be described in  simpler terms
as indicated in the following:

\begin{prop} \label{prop:example-AMM}  
Let $\Gamma$ be the AMM quasi-symplectic groupoid
$ (G \times G \toto G, \omega+\Omega)$, 
where $G$ is a connected and simply-connected Lie group equipped with 
an ad-invariant non-degenerate symmetric bilinear form. 
Let $(X\stackrel{J}{\to}G, \omega_X)$ be a quasi-Hamiltonian
$G$-space. Assume that $\omega+\Omega$ is an integral
3-cocycle in $Z^3_{dR}((G\times G)\lcom , \zz )$.
Then the pair $(\omega_X,\omega + \Omega)$
satisfies  the integrality condition
if and only if   $\forall C\in Z_2(X,{\mathbb Z})$ and 
$D \in C_3(G,{\mathbb Z})$ such that $dD=J_* (C)$,
\begin{equation}
\label{eq:AM}    
 \int_{C}\omega_X - \int_{D} \Omega \, \, \, \in \, \, \, {\mathbb Z} .
\end{equation} 
Note that such  $D$ always exists
for any $C \in Z_2(X,{\mathbb Z})$. 
\end{prop}
\begin{pf}
Note that we have the following commutating diagram
of  groupoid homomorphisms:
\begin{equation}\label{eq:Ji}\begin{array}{ccc}
   X\lcom   & \stackrel{i}{\to}& (G\times X)\lcom \\
   \downarrow J &  & \downarrow J \\
G\lcom   & \stackrel{J}{\to}&  (G\times G)\lcom,
\end{array}\end{equation}
where $X\lcom$ and $G\lcom$ are spaces $X$ and $G$ are considered
as groupoids, while $(G\times X)\lcom$ and 
$(G\times G)\lcom$ are the  transformation groupoids.
Thus one direction  is obvious.

Conversely, according to Eq.~(\ref{ref:sing-hom-grou-2}),
 we have $H_2((G \times G)_\com,{\mathbb Z})=0$.
Therefore $$\Ker(J_*)= H_2\big( (G \times M)_\com,{\mathbb Z}\big).$$
By Lemma \ref{lem:trans-cycle}, for any class
 ${\mathcal C} \in  H_2(\big(G \times X)_\com,{\mathbb Z}\big)$,
 there exists $ C\in Z_2(X,{\mathbb Z})$
such that ${\mathcal C} =i_* [C] $. Since $H_2(G,{\mathbb Z})=0$,
there always exists  $D \in C_3(G, \zz )$
such that $ J_*(C)=d D$. Hence 
$$J_* (i_* C)= i_* (J_* C) =i_* d D =\delta (i_* D ).$$
Now it is clear that
$$\int_{i_*  C}\omega_X - \int_{i_* D} (\omega +\Omega)
=\int_{C} i^*\omega_X-\int_{D} i^*(\omega +\Omega)
=\int_{C}\omega_X-\int_{D}\Omega. $$
The conclusion thus follows from  Lemma \ref{lem:simpl}.
\end{pf}

\begin{numrmk}
 Note that  Eq. (\ref{eq:AM}) coincides with the
quantization condition of  Alekseev--Meinrenken \cite{AM}. See also
\cite{Sh}. For the case of conjugacy classes, see \cite{M, Wendt}.
\end{numrmk}

As an immediate consequence of Proposition \ref{prop:example-AMM},
we have the following: 

\begin{cor}\label{cor:point}
Let $\gm$ be the  AMM quasi-symplectic
groupoid $(G \times G \toto G, \omega +\Omega )$. Then
$1 \in G$,  considered as a quasi-Hamiltonian $G$-space,
satisfies the integrality condition.
\end{cor}

 Let us consider the case of Example \ref{ex:t*G} 
where  $\Gamma$ is the symplectic groupoid  $T^*G \toto \Gg^* $.
In this case, we recover a  well-known result of Guillemin--Sternberg \cite{GS}.

 \begin{prop} \label{prop:example-integ}
Let $\Gamma$ be the symplectic groupoid  $(T^*G \toto \Gg^*, \omega )$,
 where $G$ is a connected and  simply-connected Lie group.
Let $J:X \to \Gg^*$ be a momentum map for  a 
Hamiltonian $G$-space $(X,\omega_X)$ as in Example \ref{ex:t*G}.
The pair $(\omega_X,\omega )$
satisfies the integrality condition
if and only if $\omega_X$ is an integral $2$-form.
\end{prop}
\begin{pf}
According to Eq. (\ref{ref:sing-hom-grou-1}), we have
 $H_2((T^*G)_\com,{\mathbb Z})=0$.
Therefore for any $C \in Z_2((G \times X)_\com,{\mathbb Z})$
there exists $D\in C_3((T^*G)_\com,,{\mathbb Z}) $  such that
 $J_* (C)=\delta D $.
By Lemma \ref{lem:trans-cycle}, we may  assume that
 $ C \in Z_2(X,{\mathbb Z} )$.
Since $ H_2(\Gg^*,{\mathbb Z})~=~0$, 
we can assume that $ D \in C_3(\Gg^*,{\mathbb Z})$.
Since $\Omega=0$, 
 the integrality condition of Eq.~(\ref{eq:integrable}) thus reads
    $ \int_C \omega_X \in {\mathbb Z}  $.
\end{pf}

In particular,  a  coadjoint orbit  ${\mathcal O}\subset \Gg^*$,
endowed with the Kirillov-Kostant-Souriau symplectic
structure $\omega_{\mathcal O}$,  
satisfies the integrality condition
if and only if $\omega_{\mathcal O}$ is an integral
$2$-form.


\subsection{Integrality condition and Morita equivalence}

In general, a pre-quasi-symplectic groupoid
may not be exact, as in the case of the AMM-quasi-symplectic groupoid
for instance.
In such a case, one must pass to a Morita equivalent
pre-quasi-symplectic groupoid in order to construct a 
prequantization. According to Theorem \ref{thm:morita},
 Morita equivalent   quasi-(pre)symplectic groupoids yield equivalent momentum
map theories in the sense that there is a
bijection between their (pre)-Hamiltonian $\gm$-spaces,
and the classical intertwiner spaces are  independent of
Morita equivalence \cite{Xu}.

More precisely,
given a pre-quasi-symplectic groupoid $(\gm \toto M, \omega+\Omega)$,
where $\Omega$ may not be exact, 
one can choose a  surjective submersion $N \stackrel{p}{\to} M$
 and consider the pull-back  groupoid  $\Gamma[N]\toto N$  of $\gm \toto M$ via $p$.
Then   $(\Gamma[N]\toto N, p^*\omega + p^* \Omega)$ is again a  
pre-quasi-symplectic groupoid. 
Moreover, if  $(X \stackrel{J}{\to} M, \omega_X)$ is a 
pre-Hamiltonian $\gm$-space, then $(X_N  \stackrel{J_N}{\to} N, p^* \omega_X)$ is a pre-Hamiltonian 
 $\gm [N]$-space,
where $X_N=X \times_M N $, and  $p:X_N\to X$ and $ J_N: X_N\to N$ are
 the projections to the
first and second components, respectively. The following
proposition indicates that integrality condition is  preserved
under this pull-back procedure.

\begin{lem} 
\label{cor:int1}
The pair $(\omega_X,\omega+ \Omega)$
satisfies the integrality condition 
 if and only if $\big(p^* \omega_X, p^* \omega +p^* \Omega\big)$
 satisfies the integrality condition. 
\end{lem}
\begin{pf}
By abuse of notation, we  use the same letter $p$ to denote the
groupoid homomorphisms from
$\gm[N]\times_N X_N \toto X_N$ to $\gm \times M\toto M$,
and from  $\gm[N]\toto N$ to $\gm \toto M$, both of which
are   Morita morphisms.

For any $C' \in Z_2((\gm[N]\times_N X_N)_\com,{\mathbb Z})$
and $D' \in
C_3(\gm[N]_\com,{\mathbb Z})  $ with $J_*(C' )=\delta D'$,
we have
\begin{equation}\label{eq:besoin} 
\int_{C'} p^* \omega_X -\int_{D'}p^* (\omega + \Omega) =
\int_{C}  \omega_X -\int_{D} (\omega + \Omega),\end{equation}
where   $C= p_*(C')$ and $ D=p_*(D')$  clearly satisfy  $J_*(C)=\delta D$.

 Assume that the pair $(\omega_X,\omega+ \Omega)$
satisfies the integrality condition. Then Eq.~(\ref{eq:besoin}) 
implies immediately that so too does the pair
  $(p^* \omega_X,p^*\omega+ p^*\Omega)$.

Conversely, if $\big(p^* \omega_X, p^* \omega +p^* \Omega\big)$
 satisfies the integrality condition, then $\omega+\Omega$
must be an integral cocycle.
Now we have the commutative diagram:
$$ \begin{array}{ccc} 
H_2((\gm[N]\times_N X_N )_\com,{\mathbb Z}) &\stackrel{p_*}{ \to} &
H_2((\gm \times_M X)_\com,{\mathbb Z})\\
\downarrow J_{N*} & & \downarrow J_*\\
 H_2(\gm[N]_\com,{\mathbb Z})  &  \stackrel{p_*}{\to} &   H_2(\gm_\com,{\mathbb Z}),           \end{array} $$
where  the horizontal arrows are isomorphisms. Therefore
$$ p_* :H_2((\gm[N]\times_N X_N )_\com,{\mathbb Z})\to 
H_2((\gm \times_M X)_\com,{\mathbb Z})$$
 induces an isomorphism
from  $\Ker (J_{N*})$ to $\Ker (J_*)$.
This implies that  any class in $\Ker (J_*)$
 has a representative of the form
$C=p_*(C')$ where $C'=\delta D'$ for some $D' \in 
C_3((\gm[N]\times_N X_N)_\com,{\mathbb Z})$.
Let $D = p_*(D')$. By Eq.~(\ref{eq:besoin}),
we see that
if the pair $(p^* \omega_X,p^* \omega+ p^*\Omega) $ 
satisfies the integrality condition then
  $  \int_{C}  \omega_X -\int_{D} (\omega + \Omega) \in \zz$.
By Lemma \ref{lem:simpl}, we conclude that $(\omega_X,\omega+\Omega)$
satisfies the integrality condition.
\end{pf}

\begin{cor}\label{cor:moritainteg}
Let $({G}\toto G_0 , \omega_G +\Omega_G )$
and $({H}\toto H_0 , \omega_H +\Omega_H )$ be Morita equivalent
pre-quasi-symplectic groupoids. Assume that
$(F\to G_0, \omega_F)$ and $(E\to H_0, \omega_E)$
are a pair of corresponding pre-Hamiltonian spaces.
Then $(\omega_F, \omega_G +\Omega_G )$ satisfies 
the integrality condition if and only if
$(\omega_E, \omega_H +\Omega_H )$ satisfies
the integrality condition.
\end{cor} 
\begin{pf}
It suffices to prove this assertion for a Morita morphism
of pre-quasi-symplectic groupoids. By Lemma \ref{cor:int1},
it remains to prove that the integrality
condition is   preserved by gauge transformations of
the first type, which can be easily checked.
\end{pf}

As a consequence, given a pre-quasi-symplectic groupoid
 $(\gm \toto M, \omega+\Omega)$,
where $\Omega$ may not be exact,
one can choose a  surjective submersion $N \stackrel{p}{\to} M$ such that
$p^* \Omega\in \Omega^3 (N)$ is exact
 and   replace
$(\gm \toto M, \omega+\Omega)$ by a Morita equivalent
exact pre-quasi-symplectic groupoid $(\Gamma[N]\toto N, p^* \omega+ p^*\Omega )$.
Usually, one takes $N:=\coprod U_i \to M$, where ${\mathcal{U}}=
 (U_i)$ is an open cover of   $M$. Then the pull-back
 pre-quasi-symplectic groupoid is
$(\Gamma[\mathcal{U}]\toto  \coprod U_i,
 \omega|_{\Gamma[\mathcal{U}]}+ \Omega|_{U_i})$, 
where $\Gamma[\mathcal{U}]$, as a manifold, can be identified with 
the disjoint union $\coprod \gm_{U_i}^{U_j}$.
Lemma \ref{cor:int} guarantees that  the integrality condition
always holds no matter which surjective submersion
(or open covering) $N\to M$ is taken  as long as  the initial pair
$(\omega_X,\omega+ \Omega)$
satisfies the integrality condition, and therefore
one can always construct a compatible  prequantization.

Applying the above discussion to the AMM quasi-symplectic
groupoid $(G\times G\toto G, \omega +\Omega )$ and
using  Theorem \ref{cor:quant} 
groupoid, we are led to:
 
\begin{cor}
Let $(X\stackrel{J}{\to} G, \omega_X )$ be a quasi-Hamiltonian
$G$-space. The following are equivalent
\begin{enumerate}
\item There exists a compatible prequantization
$\coprod R_{ij} \to (G\times G)[{\mathcal U}]\toto \coprod U_i $
and $\coprod L_i \to \coprod X|_{U_i} $, where
$(G\times G)[{\mathcal U}]\toto \coprod U_i $ is the
 pullback quasi-symplectic groupoid of the AMM
groupoid using any open covering
of $G$ such that $\forall i , \ \Omega|_{U_i}$ is an   exact form.
\item  The integrality condition of   Eq. (\ref{eq:AM}) holds.
\end{enumerate} 
\end{cor} 

\subsection{Strong integrality condition}
\begin{defn}
Let $(\Gamma \toto M, \omega + \Omega)$ be a
 pre-quasi-symplectic groupoid.
A pre-Hamiltonian $\gm$-space $ (X \to M, \omega_X) $
is said to satisfy the {\em strong integrality condition}  if
\begin{enumerate}
\item it satisfies the integrality condition; and
\item the map $J^* : H^2_{dR}(\gm_\com) \to H^2_{dR}((\Gamma \times_M X)_\com)$
vanishes.
\end{enumerate}
\end{defn}

The following result follows from Theorem \ref{cor:quant}.

\begin{prop}
Let $(\Gamma \toto M,\omega+\Omega)$
be  an   exact proper  pre-quasi-symplectic groupoid,
and  $(X \stackrel{J}{\to} M, \omega_X)$  a
pre-Hamiltonian $\gm$-space.
Then $(X \stackrel{J}{\to} M, \omega_X)$ satisfies
the strong integrality condition  if and only if
for any prequantization of $(\Gamma \toto M,\omega+\Omega)$,
$X$ admits a compatible prequantization.
\end{prop} 
\begin{pf}
If $(X \stackrel{J}{\to} M, \omega_X)$ satisfies
the strong integrality condition, it is clear from  Theorem \ref{cor:quant}
that $X$ admits a compatible prequantization for any
 prequantization of $(\Gamma \toto M,\omega+\Omega)$.

Conversely, given any prequantization $(R\to \gm \toto M, 
\theta +B)$, 
 $J^*(\theta + B)-\omega_X$ must be an integral
2-cocycle in $ Z^2_{dR}((R \times_M X)_\com,{\mathbb Z})$. 
Note that if  $\theta + B $ is a pseudo-connection, so is
$\theta + B + \pi^* \Xi, \ \forall \Xi \in Z^2_{dR}(\gm\lcom) $.
Since the subset of integral classes  $ Z^2_{dR}(R \times_M
 X_\com,{\mathbb Z})$
is   discrete, then
  $J^*(\theta + B)-\omega_X + J^* \pi^* \Xi  $ being
 an integral  cocycle for all $\Xi  $
 implies that 
$[J^* \smalcirc \pi^* (\Xi)]=0$. 
In other words,  the map
$ J^* \smalcirc \pi^*: H^2_{dR}(\gm\lcom ) \to H^2_{dR}((R \times_M X)\lcom) $
is the zero map. Since $R \times_M X \to \Gamma \times_M X \toto X $
defines a trivial gerbe according to Proposition \ref{pro:2},
the map $\pi^*: H_{dR}^2(\Gamma \times_M X) \to H_{dR}^2 (R \times_M
X)  $ is injective. From the identity $J^* \smalcirc \pi^*= \pi^* \smalcirc J^* $
 and the fact that $\pi^*$ is injective,
it follows that $ J^*:  H^2_{dR}(\gm\lcom) \to 
H^2_{dR}((\gm \times_M X)\lcom)$ must vanish. \end{pf}

The following proposition is an
 analogue of Corollary \ref{cor:moritainteg}.

\begin{prop}\label{prop:moritainteg2}
Let $({G}\toto G_0 , \omega_G +\Omega_G )$
and $({H}\toto H_0 , \omega_H +\Omega_H )$ be Morita equivalent
pre-quasi-symplectic groupoids. Assume that
$(F\to G_0, \omega_F)$ and $(E\to H_0, \omega_E)$
are a pair of corresponding pre-Hamiltonian spaces.
Then $(\omega_F, \omega_G +\Omega_G )$ satisfies 
the strong integrality condition if and only if
$(\omega_E, \omega_H +\Omega_H )$ satisfies
the strong integrality condition.
\end{prop} 
\begin{pf} By Corollary 
\ref{cor:moritainteg}, we just have to check that
Condition (2) in the definition of strong integrality
condition is invariant under Morita equivalence.
This follows immediately
from the commutativity of the   diagram
$$\begin{array}{ccc}
  H^2_{dR}(G_\com) &    \simeq   &   H^2_{dR}(H_\com)  \\
 \downarrow &        & \downarrow\\ 
 H^2_{dR}((G \times_{G_0} F)_\com)  &\simeq  &
 H^2_{dR}((H \times_{H_0} E)_\com),     \\ 
\end{array}$$
where the horizontal arrows  are the natural isomorphism
between the de Rham cohomologies of two Morita
equivalent groupoids.
\end{pf}

\begin{numrmk}\label{ex:h1triviakl}
\begin{enumerate}
\item
If the groupoid $ \Gamma $ satisfies $H^2_{dR}(\Gamma_\com,{\mathbb Z})=0$,
then Condition (2) in the definition
of the strong integrality is satisfied for any pre-Hamiltonian
 $\gm$-space.
In this case, a   pre-Hamiltonian space
satisfies the integrality condition
if and only if it satisfies the strong integrality condition.
\item If $G$ is a connected and
 simply-connected  Lie group, then $H_2((G\times G)_\com,{\mathbb Z})=0$.
Therefore, any quasi-Hamiltonian $G$-space satisfying
the integrality condition
must satisfy the strong integrality condition.
\end{enumerate}
\end{numrmk}

The following proposition summarizes the results
of this section.

\begin{prop} 
Let $({\gm}\toto{M}, \omega +\Omega )$ be an exact,
proper,  pre-quasi-symplectic   groupoid,
and   $(X_k\stackrel{J_k}{\to} M, \omega_k), \ k=1, 2,$ be 
 pre-Hamiltonian $\gm$-spaces.
Assume that $(X_1\stackrel{J_1 }{\to} M, \omega_1)$ satisfies
 the integrality condition while $(X_2\stackrel{J_2}{\to} M, \omega_2)$
 satisfies the strong integrality condition. Then
there exists   a prequantization of $({\gm}\toto{M}, \omega +\Omega )$
and compatible prequantizations  of both $X_1$ and $X_2$. Therefore
the classical intertwiner space $\overline{X_2}\times_\gm X_1$
is quantizable.
\end{prop}

Applying this result to the case of the AMM quasi-symplectic
groupoid, we have the following

\begin{cor}
Let $G$ be  a connected and simply-connected compact
 Lie group equipped with
an ad-invariant non-degenerate symmetric bilinear form, and
 $(X\stackrel{J}{\to}G, \omega_X)$  a quasi-Hamiltonian
$G$-space. Assume that $\omega_X$ satisfies the
integrality condition as in Eq. (\ref{eq:AM}).
Then the reduced symplectic manifold $J^{-1}(1)/G$ is prequantizable,
and the prequantization can be constructed using
the prequantization of the   AMM quasi-symplectic groupoid
$(G\times G)[{\mathcal U}]\toto \coprod U_i $ (more precisely
 the pull-back groupoid of the AMM quasi-symplectic groupoid) together
with a compatible prequantization of the Hamiltonian space
$(\coprod X|_{U_i}\to  \coprod U_i, \omega_X|_{U_i})$, 
where ${\mathcal U}= (U_i)_{i\in I}$ is some open covering
of $G$ such that $ \Omega|_{U_i}$ is  exact $\forall i \in I$.
\end{cor}

\section{Appendix}

We denote  then by $C_{S^1}$  the canonical cycle
in $ C_1(S^1,{\mathbb Z})$ that generates $H_1(S^1,{\mathbb Z})=\zz $.
If we consider $C_{S^1} $ as an element of $C_2(S^1_\com,{\mathbb Z}) $,
$[C_{S^1}]$ generates $ H_2(S^1_\com,{\mathbb Z})\simeq {\mathbb Z}$.
For  any point $p $ in a manifold   $N$, we denote by $C_p$ 
the constant map from $S^1$ to  $\{p\}$ and consider it
 as an element of $C_1(N,{\mathbb Z})$.
Assume that  $R \to \gm \toto M$ is 
an $S^1$-central extension of groupoids. For any $m \in M$,
 let 
\begin{equation} \label{defZm} 
Z_m = f_{m*}(C_{S^1})\in C_2(R_\com,{\mathbb Z}),
\end{equation}
where $f_m : S^1 \to R$ is defined by Eq.  (\ref{eq:plonge}).
More generally, for any $r \in R$, let
 $f_r:S^1 \to R$ be the map $\lambda \to \lambda\cdot r$,
 and set  
$$Z_r=f_{r*} (C_{S^1})-C_r \in C_2(R_\com,{\mathbb Z}). $$

\begin{prop}\label{lem:precisions} Let $R \to \gm \toto M$ be
an $S^1$-central extension. Assume that $M/\gm$ is connected. 
\begin{enumerate}
\item The class
 $[Z_m] \in H_2(R_\com,{\mathbb Z})$ does not depend 
on the choice of $m \in M$. Because of this, 
 we will drop the subscript $m$ and denote
this class   simply by $[Z]$.
\item For any $r \in R$, $Z_r$ is a cycle and $[Z_r]=[Z]$;
\item  The natural map $\pi_*:H_2(R_\com,{\mathbb Z}) \to 
H_2(\Gamma_\com,{\mathbb Z})$ is surjective.
\item Its kernel $\Ker(\pi_*)$ is generated by
 $[Z]$. 
\item The following  identity holds:
 $$ \int_{Z_m} \theta =1 . $$
\end{enumerate}
\end{prop}

Before we prove this proposition, we first need a lemma.
Given any point $p\in N$,
 we will denote by $C_{p^{(k)}}$ the chain in $C_k(N,{\mathbb Z})$
defined by the constant path $\Delta_k \to  \{p\}$.

\begin{lem} Let $R \stackrel{\pi}{\to} \Gamma \toto M$ be an $S^1$-central extension.
\label{lem:details}
\begin{enumerate}
\item  Any element $E $ in   $C_0(R,{\mathbb Z})$
with $\pi_*(E)=0$ can be written of  the form $E=\delta D'$, 
where $D' \in C_1(R,{\mathbb Z}) $ satisfies $\pi_* (D')=0$.
\item $\pi_*:C_\com(R_\com,{\mathbb Z}) \to C_\com(\gm_\com,{\mathbb Z})$
is a surjective map.
\item Any element in  the kernel of  $\pi_*: H_k(R_\com,{\mathbb Z}) \to
H_k(\gm_\com,{\mathbb Z})$ has a representative $C \in Z_k(R_\com,{\mathbb Z})$ with $\pi_*(C)=0$.
\item Any element $C $ in   $C_0(R_2,{\mathbb Z})$
with $\pi_*(C)=0$  is of the form $C =d  D'$, 
where $D' \in C_1(R_2,{\mathbb Z}) $ satisfies $\pi_*(D')=0$.
\item For any cycle $C' \in C_1(R,{\mathbb Z})$ such that  $\delta C'=0$
and  $\pi_*(C')=0$, we have
 $$ [C']= \sum_{i\in I} k_i [Z_{r_i}]$$
for some finite set $I$, $k_i \in \zz$  and $r_i \in R$.
\end{enumerate}
\end{lem}
\begin{pf}
1) The kernel of  $\pi_*:C_0(R,{\mathbb Z})\to C_0(\gm,{\mathbb Z}) $
is  generated by elements of the form $p -q$,
where $p$ and $q$ are two points in
 the same fibre of $R \stackrel{\pi}{\to} \gm$. Hence,
 it suffices to prove the claim for such a generator.  

Let $D:\Delta_1 \to R$ be a path 
in the fiber $\pi^{-1}(p)$ satisfying $d D= p -q $.
Set  $D'=D- C_{p^{(1)}} $.
Clearly, the identities $d D'=p -q$   and $\pi_*(D')=0$ hold.
Moreover, from $\pi_*(D')=0$, it  follows that
  $\partial  D'= s_* \pi_*(D')-t_* \pi_* (D')$. Hence $\delta D'=p-q$

2) Since  the projections $\pi:R_k \to \Gamma_k$ are surjective
submersions
with fibers isomorphic to $k$-dimensional torus, 
 all the maps $\pi_*:C_l(R_k,{\mathbb Z}) \to C_l(\gm_k,{\mathbb Z}) $
are onto for all $k,l \in {\mathbb N}$.

3) Let $ C' \in Z_k(R_\com,{\mathbb Z})$ be a cycle with $\pi_*[C']=0$.
By definition, there exists $D \in C_{k+1}(\gm_\com,{\mathbb Z})$
such that  $\delta D =\pi_*(C') $. 
By 2),  there exists $D'$ in  $C_{k+1}(R_\com,{\mathbb Z})$ 
such that  $\pi_*(D')=D$. Set $C:= C' -\delta D'$.
We have $[C]=[C'] $ and $\pi_*(C)= 0$.

4) The kernel of $\pi_*: C_0(R_2,{\mathbb Z})\to C_0(\gm_2,{\mathbb Z})$
is  generated by elements of the form $p-q $,
where $p$ and $q$ are two points on  the same fiber of 
$R_2 \to \gm_2$. It thus suffices to show this property for such  generators.

Let $D:\Delta_1 \to R_2$ be a path in the fiber over  $\pi(p)$
 such that $d D =p -q$.
Let $D'=D- C_{p^{(1)}}$. Thus
$$\pi_*(D')=\pi_*(D- C_{p^{(1)}})=C_{\pi (p)^{(1)}}-C_{\pi(p)^{(1)}}=0,\ \ \ \
\mbox{ and } \ \ 
  d D'=p-q .$$

5)  For simplicity, we call  those chains
in $C_1(R,{\mathbb Z}) $ of the form $ C -C_{p^{(1)}}$ {\em fibered $1$-chains},
where $p\in R$ is a point
and $C : \Delta_1 \to \pi^{-1}(\pi(p))$
is a path  in the fiber through
the   point $p$.  Any fibered 1-chain 
is in the kernel of $\pi_*$ and hence lies in the kernel
of $ \partial$. If  a fibered 1-chain $E$ in a given fiber 
satisfies $dE=0$,
then $[E]=k [Z_r]$ for some $r \in R$  and  $k \in {\mathbb Z}$.   
As a consequence, if a linear combination 
of fibered 1-chains $F$ is a cycle in $C_1(R, {\mathbb Z})$,
 then it is clear that $[F]=\sum_{i\in I} k_i [Z_{r_i}]$
for some finite set $I$, $k_i\in \zz$ and $r_i \in R$.

Now the kernel of $\pi_*: C_1(R,{\mathbb Z})\to C_1(\gm,{\mathbb Z})$
is generated by elements of the form 
$C_0 -C_1 $,  where $C_i,i=0,1,$ are  paths $\Delta_1 \to R$
satisfying  $\pi_* (C_0)=\pi_*(C_1)$. 
Thus there is a map  $\gamma: \Delta_1 \to S^1 $
such that  $ C_0(t) = \gamma(t) \cdot C_1(t)$   $\forall t \in \Delta_1$.
Let $\tilde{\gamma}: [0,1]\times \Delta_1 \to S^1$
be a map with $\tilde{\gamma}(0,t)=1 $ and $\tilde{\gamma}(1,t)=\gamma(t) $.
Let us define 
two maps  $D_1$ and $\hat{D}: [0,1]\times \Delta_1 \to R$ by 
$ D_1(s,t)= \tilde{\gamma}(s,t) \cdot C_0(t)$ and   $\hat{D}(s,t)=C_0(t)$.
Set  $D:= D_1- \hat{D} $.
  We have $\pi_* ( D)=0$ and therefore $\partial D=0$.
Moreover, by construction, $C_0 -C_1 +\delta D=C_0 -C_1 +d D $ is
 the sum of two 1-fibered chains:
one in the fiber through $C_0(0)$ and another in the fiber through $C_0(1)$.
The conclusion thus follows.
\end{pf}

\begin{pff}
1) and  2). It is clear that  if $m$ and $n$ are in the 
same connected component of $M$, then
$[Z_m] =[Z_n]$.
Now   by the  definition of $Z_r$, we have
\be
 &&d Z_r =d\big(f_{r*} (C_{S^1})\big)- d C_r=0, \\
 && s_*(Z_r)=s_*\big(f_{r*}(C_{S^1})-C_r \big)=C_{s(r)}-C_{s(r)} =0, \ \mbox{and}\\
&& t_*(Z_r)=t_*\big(f_{r*}(C_{S^1})-C_r \big)=C_{t(r)}-C_{t(r)} =0.
\ee
Therefore $\delta (Z_r)=0$.  Consider
 the map  $D:S^1 \to R_2$ defined by
 $ \lambda \mapsto \big(f_r(\lambda), f_{t(r)}(\lambda^{-1})\big) $.
 We have $d D =0$ and
 $ \partial D = f_{r*} (C_{S^1})  - C_r  - f_{t(r)*} (C_{S^1})$. 
Hence we have 
$ [Z_r]=[Z_{t(r)}], \ \forall r \in R$.
Similarly, we have $[Z_r]=[Z_{s(r)}], \ \forall r \in R$.
Since $M/\Gamma$ is connected,  1) and 2) follow.

3)
Let  $C \in Z_2(\Gamma_\com,{\mathbb Z})$ 
be any 2-cycle.
According to Lemma \ref{lem:details} (2),
there exists $D\in C_2(R_\com,{\mathbb Z})$
with $\pi_*(D)=C $. 

In general, $\delta D \neq 0 $. However since the restriction 
of $\pi$ to $M$
is the identity map, we have $\partial~D_1~-~d~D_2~=~0$ 
and thus $\delta D =\partial D_0-dD_1$, where $D=D_0 +D_1 +D_2, \ 
D_i \in C_i (R_{2-i}, \zz)$.
Therefore $ \delta D$ is an element of $C_0(R,{\mathbb Z})$
and  $\pi_* (\delta D)=\delta \pi_* (D)=\delta C=0$.
By Lemma \ref{lem:details} (1), 
there exists
 $ D' \in C_2(R_\com,{\mathbb Z})$
with $\pi_*(D')=0$ and $\delta D'=\delta D $.
Therefore it follows that $D-D' $ is a cycle in
 $Z_2(R_\com,{\mathbb Z})$ and
$$\pi_*([D-D'])=[\pi_*(D)]- [\pi_*(D')] =[C]-[0]=[C].$$

4)  According to  Lemma \ref{lem:details} (3),
 any class in $\Ker(\pi_*)$
has a representative $C$ such that  $\pi_*(C)=0$ and therefore
is of   the form $C_0+C_1$, where $C_0 \in C_0(R_2,{\mathbb Z})
$ and $C_1
\in   C_1(R,{\mathbb Z})$ 
satisfy $\pi_*(C_0)=0 $ and $ \pi_*(C_1)=0$.
According to  Lemma \ref{lem:details} (4), 
there exists $D'\in C_1(R_2,{\mathbb Z})$
 with $\pi_* ( D' )=0 $ such that $C_0 = dD'$.
Consider now $C'=C-\delta D' \in  C_1(R,{\mathbb Z})$. We have 
 $$ \delta C' = \delta C-\delta^2 D'=0,  \ \ 
[C']=[C] , \ \ \ 
\pi_*(C')=0.  $$

According to  Lemma \ref{lem:details} (5), we have

\begin{equation}\label{eq:suppl5} [C']=\sum_{i\in I} k_i [Z_{r_i}] \end{equation}
for some 
  finite set $I$,
 $k_i \in {\mathbb Z}$ and   $r_i \in R  $.
 From Eq.~(\ref{eq:suppl5}), it  follows that
$[C]=[C']= \sum_{i\in I} k_i [Z_{r_i}]$ .
By  Lemma \ref{lem:precisions} (2),  we have
  $ [C]=(\sum_{i\in I} k_i) [Z]$.

5) holds because $\theta$ is a connection $1$-form of the $S^1$-principal bundle
   $R\to \Gamma$. \end{pff}

\end{document}